\documentclass[11pt,leqno]{amsart}
\usepackage[english]{babel}  

\usepackage{geometry}
\usepackage{graphicx}     
\usepackage{xcolor}       
\usepackage{pstricks-add} 
\usepackage{enumitem}     
\usepackage{stmaryrd}     
\SetSymbolFont{stmry}{bold}{U}{stmry}{m}{n}
\usepackage{morefloats}

\usepackage{amsmath, amsthm}
\usepackage{amsfonts, amssymb}
\usepackage{bm}               
\usepackage{hyperref}

	\newtheorem*{theorem*}{Theorem}

	\newtheorem*{proposition*}{Proposition}

\theoremstyle{remark}
\newtheorem{remark}{\textbf{Remark}}[]

\newtheorem*{example}{Example}

\theoremstyle{definition}

\DeclareMathOperator{\vect}{Span}

\DeclareMathOperator{\tr}{Trace}

\DeclareMathOperator{\diverg}{div}

\newcommand{\ve}{\varepsilon}
\newcommand{\eg}{\emph{e.g. }}
\newcommand{\ie}{\emph{i.e. }}
\newcommand{\RR}{\mathbb{R}}

\newcommand{\TT}{\mathbb{T}}

\newcommand{\Q}{\mathcal{Q}}
\newcommand{\G}{\mathcal{G}}
\newcommand{\Tau}{\mathcal{T}}
\newcommand{\LL}{\mathcal{L}}

\title{A Hierarchy of Hybrid Numerical Methods for Multi-Scale Kinetic Equations}

\author[F. Filbet]{Francis Filbet}
\address{Francis Filbet \\
Universit\'e de Lyon \\
CNRS UMR 5208 \\
Universit\'e Lyon 1 \\
Institut Camille Jordan \\
43 blvd. du 11 novembre 1918 \\
F-69622 Villeurbanne cedex \\
France.}
\email{filbet@math.univ-lyon1.fr}

\author[Th. Rey]{Thomas Rey}
\address{Thomas Rey \\
Center of Scientific Computation and Mathematical Modeling (CSCAMM) \\
The University of Maryland \\
College Park, MD, 20742-4015\\
USA}
\email{trey@cscamm.umd.edu}

	

\begin{document}
  
	\begin{abstract}
  	In this paper, we construct a hierarchy of hybrid numerical methods for multi-scale kinetic equations based on moment realizability matrices, a concept  introduced by Levermore, Morokoff and Nadiga in \cite{LevermoreMorokoffNadiga:98}. Following such a criterion,  one can consider hybrid scheme where the hydrodynamic part is given either by the compressible Euler or Navier-Stokes  equations, or even with more general models, such as the Burnett or super-Burnett systems.
	
	\end{abstract}
  \maketitle

  \tableofcontents

  \section{Introduction}\label{sec1}
    \setcounter{equation}{0}
	  Many engineering problems involve fluids in transitional regimes (micro-electro-mechanical systems, space shuttle reentry, ...). In these cases, the Euler or Navier-Stokes-like fluid description breaks down, typically due to shocks or boundary layers, and the use of a kinetic model is needed to describe accurately the system. 
	  Nevertheless, this type of mathematical description is computationally expansive to simulate, and it is desirable to use it only locally in space.
	  The goal of this paper is then to design a so-called \emph{hybrid kinetic/fluid schemes} with an automatic domain-decomposition criterion allowing to identify accurately the fluid and kinetic zones.
	  For the sake of computational efficiency, we will give a decomposition which minimizes the size of the kinetic layer, allowing to take advantage of the low computational cost of numerical methods for  fluid systems. 
	  As far as possible, this method will also be \emph{non-intrusive} for the solvers. More precisely, except for the implementation of the domain decomposition indicators, it will be independent on the kinetic and fluid solvers, which won't necessitate deep modifications.

    
	  We are interested in this article in collisional gases, and then we shall consider Boltzmann-like collisional kinetic equations. More precisely, for a given nonnegative initial condition $f_0$, we will study a particle distribution function $f^\ve = f^\ve(t,x,v)$, for $t \geq 0$, $x \in \Omega \subset \mathbb{R}^{d_x}$ and $v \in \mathbb{R}^{3}$, solution to the initial-boundary value problem
		\begin{equation} 
			\label{eqCollision}
			\left\{ \begin{aligned}
			  & \frac{\partial f^\ve}{\partial t} + v \cdot \nabla_x f^\ve \,=\, \frac{1}{\ve}\;\Q(f^\ve), 
			  \\
		  	&  \;
			  \\
			  & f^\ve(0, x, v) = f_{0}(x,v),
			\end{aligned} \right.
		\end{equation}
		where the collision operator $\Q$ is a  Boltzmann-like operator. 
		The open  set $\Omega$ is a bounded Lipschitz-continuous domain of $\RR^{d_x}$, which means  that the model (\ref{eqCollision}) has to be supplemented with boundary conditions described later.
		
		We assume that the collision operator fulfills the three following assumptions
		\begin{enumerate}[label=\textbf{(H$\bm{_\arabic{*}}$)}, ref=\textbf{(H$\bm{_\arabic{*}}$)}]
		
      \item \label{hypConservations} 
        It preserves mass, momentum and kinetic energy:
			  \begin{equation*}
			    \int_{\RR^{3}} \Q(f)(v) \, dv = 0, \quad  \int_{\RR^{3}} \Q(f)(v) \, v \, dv = 0, \quad \int_{\RR^{3}} \Q(f)(v) \, |v|^2 \, dv = 0;
			  \end{equation*}
			 
			\item \label{hypEntropy}
			  It dissipates the Boltzmann entropy (H-theorem):
			  \begin{equation*}
					\int_{\RR^{3}} \Q(f)(v) \, \log(f)(v) \, dv \, \leq \, 0;
				\end{equation*}
				
	    \item \label{hypEquilib} 
	      Its equilibria are given by Maxwellian distributions:
	      \begin{equation*}
	        \Q(f) \, = \, 0 \quad \Leftrightarrow \quad f = \mathcal M_{\rho, \bm u, T} := \frac{\rho}{(2 \pi T)^{3/2}} \exp \left ( - \frac{|v-\bm u|^2}{2 T} \right ),
	      \end{equation*}
	      where the \emph{density}, \emph{velocity} and	\emph{temperature} of the gas $\rho$, $\bm u$ and $T$  are computed from the distribution function $f$ as  
			\begin{equation*}
	      \rho = \int_{v\in{\RR}^{3}}f(v)\,dv, \quad \bm{u} = \frac{1}{\rho}\int_{v\in{\RR}^{3}}v f(v) \, dv, 
	      \quad T = \frac{1}{3\rho} \int_{v\in{\RR}^{3}}\vert \bm{u} - v \vert^2 f(v) \,dv.
	    \end{equation*}
	      
	  \end{enumerate}  
	  Equation \eqref{eqCollision} with assumptions \ref{hypConservations}-\ref{hypEntropy}-\ref{hypEquilib}  describes numerous models such as the Boltzmann equation for elastic collisions or Fokker-Planck-Landau type equations. 
	   
		The parameter $\ve > 0$ is the dimensionless Knudsen number, that is the ratio between the mean free path of particles before a collision and the length scale of observation. It measures the rarefaction of the gas: the gas is in \emph{rarefied} or kinetic regime if $\ve \sim 1$ and in \emph{dense} or fluid regime if $\ve \ll 1$.
		Moreover, according to assumptions \ref{hypEntropy}-\ref{hypEquilib}, when $\ve \to 0$, the distribution $f^\ve$ converges (at least formally) to a Maxwellian distribution, whose moments are solution to the compressible Euler system
			\begin{equation}
				\label{eqHydroClosedEuler}
				\left\{ \begin{aligned}
				  & \partial_t \rho  + \diverg_x  (\rho \, \bm{u} ) \,=\, 0, 
				  \\
				  \,
				  \\
				  & \partial_t(\rho \, \bm{u} ) + \diverg_x  \left(\rho \, \bm{u} \otimes \bm{u} \,+\, \rho \, T  \,{\rm\bf I}\right) \, =\, \bm{0}_{\RR^{3}}, 
				  \\
				  \,
				  \\
				  & \partial_t E + \diverg_x \left ( \bm{u} \left ( E +\rho \, T\right )  \right ) \,=\, 0.
				\end{aligned} \right.
			\end{equation}
	  This type of limit provides a ``contraction of the kinetic description'' \cite{ellis:1975}, the so-called \emph{hydrodynamic limit}, which is at the basis of the hybrid methods.

    Let us give some examples of operators which are subject to these hypothesis.
    
    \subsection{The Boltzmann Operator}
	    \label{subBoltzOp}
	    
		  The Boltzmann equation describes the behavior of a dilute gas of particles when the only interactions taken into account are binary	elastic collisions

			\begin{align}
			  \label{eqQBoltz}
			  \Q_\mathcal{B} (f,f)(v) & \,=\, \int_{\RR^{3}} \int_{\mathbb{S}^{2}}  B(|v-v_*|,\cos \theta) \, \left[ f'_* f' \,-\, f_* f \right] \, d\sigma \, dv_*,
				\end{align}
			where we used the shorthand $f = f(v)$, $f_* = f(v_*)$, $f ^{'} = f(v')$, $f_* ^{'} = f(v_* ^{'})$.   
			The velocities of the colliding pairs $(v,v_*)$ and $(v',v'_*)$ are related by
			\begin{equation*}
			  v' \,=\, \frac{v+v_*}{2} \,+\, \frac{|v-v_*|}{2} \,\sigma, \qquad v'_* \,=\, \frac{v+v^*}{2} \,-\, \frac{|v-v_*|}{2} \,\sigma.
			\end{equation*}
			The \emph{collision kernel} $B$ is a non-negative function which by physical arguments of invariance only depends on $|v-v_*|$ and	$\cos \theta = {\widehat u} \cdot \sigma$, where ${\widehat u} =(v-v_*)/|v-v_*|$. 
									
			Boltzmann's collision operator has the fundamental properties \ref{hypConservations}, \ref{hypEntropy} and \ref{hypEquilib}.

		\subsection{The BGK Operator}
		  \label{subBGK}
		
		  Another well known collision operator which has the properties \ref{hypConservations}-\ref{hypEntropy}-\ref{hypEquilib} is the BGK operator \cite{Bhatnagar:1954}, and its Ellipsoidal Statistical (ES-BGK) extension \cite{AndriesPerthame:2000}.
		  It consists in replacing the \emph{bilinear} collision operator $\Q_\mathcal{B}$ by a \emph{nonlinear} relaxation operator, which match the same hydrodynamic limit than the Boltzmann operator.

      To this aim we  first define some macroscopic quantities of the particle distribution function $f$ such as the opposite of the \emph{stress tensor}
			\begin{equation*}
			  \Theta_f(t,x)\, =\, \frac{1}{\rho_f} \int_{{\RR}^{3}}  (v-\bm{u}_f)\otimes (v-\bm{u}_f)\,f(t,x,v) \,dv.
			\end{equation*} 
			Therefore the \emph{translational temperature} is related to the stress tensor as $T_f = \tr(\Theta_f)/3$. We finally introduce the corrected tensor
			\begin{equation*}
			  \Tau_f(t,x) \,\,=\,\, \left[(1-\beta) \, T_f \,{\rm\bf I} \,\,+\,\,\beta \,\Theta_f\right](t,x),  
			\end{equation*} 
			which can be viewed as a linear combination of the initial stress tensor $\Theta_f$ and of the isotropic stress tensor $T_f \,{\rm \bf I}$ developed by a Maxwellian distribution. 
			The parameter $-\infty< \beta < 1$ is used to modify the value of the Prandtl number through the formula 
			\[
			0 \,\leq\, {\rm Pr}  \,=\, \frac{1}{1-\beta} \,\leq\, +\infty\quad {\rm for } \quad \beta\in (-\infty\,,\,1).
			\]
			The correct Prandtl number for a monoatomic gas of hard spheres is equal to $2/3$, namely obtained here for $\beta = -1/2$ whereas the classical BGK operator, obtained for $\beta = 0$, has a Prandtl number equal to $1$.
			
			To define the ESBGK operator, we introduce a corrected Gaussian $\G[f]$ defined by
			\[
			\G[f]= \frac{\rho_f}{\sqrt{{\rm det}(2\pi\,\Tau_f)}}\,\exp\left(-\frac{(v-\bm{u}_f)\,\Tau_f^{-1}\,(v-\bm{u}_f)}{2}\right)
			\] 
			and the corresponding collision operator is now
			\begin{equation} 
			  \label{eqOpESBGK}
			  \Q_\mathcal{BGK}(f) =  \nu(\rho_f, T_f) \left( \G[f] \,-\, f \right),
			\end{equation}
			where $\nu$ is the collision frequency from the Boltzmann operator. It can be shown \cite{Struchtrup:2005} that it depends only on the kinetic density $\rho_f$ and temperature $ T_f $.

  \section{Regime Indicators}
    \label{secRegimeIndic}
    
    There are several different works about hybrid methods in the literature, the largest part relying on the same domain decomposition technique, introduced by Boyd, Chen and Chandler in \cite{Boyd:1995}.
	  This paper uses a macroscopic criterion to pass from the hydrodynamic description (easy to compute numerically, but inaccurate near shocks or boundary layers) to the kinetic one (computationally expansive but accurate in most of the situations).
	  This criterion is based on the local Knudsen number of the problem: when this quantity is below a (problem-dependent) threshold, the kinetic description is used.
	  The first practical use of this method is due to Kolobov, Arslanbekov, Aristov \textit{et al.} \cite{Kolobov:2007}, by using a discrete velocity model of the Boltzmann equation for the kinetic part, and a kinetic scheme for the hydrodynamic part.
	  It has been more recently used by Degond and Dimarco in \cite{DegondDimarco:2012}, with a Monte-Carlo solver for solving the kinetic equation and a finite volume method for the macroscopic ones.
	  
		The hydrodynamic breakdown indicator introduced by Tiwari in \cite{Tiwari:98} is also very close to the previous criterion, as it is based on the viscous and heat fluxes of the Navier-Stokes equation, through a Grad's 13-moments expansion. It has been recently used for deterministic solver by Degond, Dimarco and Mieussens in \cite{DegondDimarcoMieussens:2010} and by Tiwari, Klar et Hardt in \cite{TiwariKlarHardt:2009,TiwariKlarHardt:2012}.
		Alaia and Puppo also used it with an hybrid deterministic/kinetic solver in \cite{AlaiaPuppo:2012}.
		Finally, Dimarco, Mieussens and Rispoli used in \cite{DimarcoMieussensRispoli:2013} the same ideas to solve the more complex Vlasov-Poisson-BGK system.
		Another different approach introduced by Dimarco and Pareschi in \cite{DimarcoPareschi:2008} consists in decomposing the particle distribution in a ``central'' part containing most of the information, solved using a BGK model, and a Monte Carlo part for the tail, being able to attain any velocity.
		
		To our opinion, the main outcome of these criteria is that they are based on the macroscopic description of the model, and more precisely on the magnitude of the first spatial derivatives of the local density, temperature, heat flux and stress tensor. In consequence, they can be wrong in a situation where the fluid is far from the thermal equilibrium. 
		Let us consider for example the case of a gas which is distributed in the velocity space as a sum of two Gaussian with non-zero mean, and which is constant in space:
		\begin{equation}
		  \label{eqSumMaxwell}
	    f(x,v) = \frac{1}{2}\left ( \mathcal{M}_{1, \bm{u}_0, 1}(v) + \mathcal{M}_{1, -\bm{u}_0, 1}(v) \right ), \quad \forall x \in \TT, \, v \in \RR^3,
		\end{equation}
		for \eg $ \bm{u}_0 = (1,0,1)$. This distribution, although constant in space, is very far from the \emph{thermal equilibrium} given by Maxwellian distributions according to assumption \ref{hypEquilib}. Nevertheless, both criteria from  \cite{Boyd:1995} or \cite{Tiwari:98} would be equal to $0$, since they are based on the spatial derivatives of the hydrodynamic field, hence detecting an hydrodynamic setting. 
		
		We then need criteria in both regimes: one to know when the hydrodynamic description breaks down, and another one to know when the kinetic description is actually in hydrodynamic regime.
		Let us start by introducing some mathematical tools needed for our approach, namely the Chapmann-Enskog expansion.
	
    \subsection{The Chapmann-Enskog Expansion}			
      \label{subChapmannEnskog}
	 	
	 	  The criterion we shall use in this article was introduced by Levermore, Morokoff and Nadiga in \cite{LevermoreMorokoffNadiga:98}. It has already been used by Tiwari in \cite{Tiwari:2000} for an hybrid Euler-Boltzmann method using particle methods.
	 	  It has the main interest to depend on the closure made for obtaining the hydrodynamic model. Hence, one can consider hybrid scheme where the hydrodynamic part is given either by the compressible Euler equations, or by compressible Navier-Stokes, or even with more general models, such as the Burnett or super-Burnett systems,and we shall take advantage of this to design a hierarchy of models.
		  
		  Let us consider a solution $f^\ve$ of the collisional kinetic equation \eqref{eqCollision}. Without any closure, according to the conservative properties \ref{hypConservations} of the collision operator $\Q$, we have
		  \begin{equation}
		    \label{eqHydroGeneral}
				\left\{ \begin{aligned}
					& \partial_t \rho^\ve  + \diverg_x (\rho^\ve \, \bm{u}^\ve ) \,=\, 0, 
					\\
					& \,
					\\
					& \partial_t(\rho^\ve \, \bm{u}^\ve ) + \diverg_x \left (\int_{\RR^{3}}  v \otimes v \, f^\ve(v) \, dv\right ) \,=\, \bm{0}_{\RR^{3}}, 
					\\
					& \,
					\\
					& \partial_t E^\ve + \diverg_x \left (\int_{\RR^{3}} \frac{1}{2}| v|^2 v \, f^\ve(v) \, dv \right ) \,=\, 0,
				\end{aligned} \right.
			\end{equation}
			where we defined
			\[
			  (\rho^\ve,\bm{u}^\ve,E^\ve) \, = \, \int_{\RR^{3}} f^\ve(v) \, \left(1, v, \frac{|v|^2}{2}\right) \, dv; \quad T^\ve \, = \, \frac{1}{3 \, \rho^\ve }\int_{\RR^{3}} f^\ve(v) \, |v-\bm{u}^\ve|^2 \, dv.
			\] 
			Now, assuming that the distribution $f^\ve$ is close to equilibrium thanks to the relaxation property \ref{hypEntropy}, we can do formally the \emph{Chapman-Enskog} expansion
			\begin{equation}
			  \label{devChapEnsk}
			  f^\ve \, = \, \mathcal M_{\rho^\ve,\bm{u}^\ve,T^\ve}\left [1 + \ve \, g^{(1)} + \ve^2 \, g^{(2)} + \ldots \right ],
			\end{equation}
			where the fluctuations  $g^{(i)}$ for $i \geq 0$ designate a function that depends smoothly on the moment vector $(\rho^\ve,\bm{u}^\ve,T^\ve)^\intercal $  and any finite number of its derivatives with respect to the $x$-variable at the same point $(t, x)$, and on the $v$-variable. According to \ref{hypConservations}, it verifies
			\[
			  \int_{\RR^{3}} g^{(i)}(v) \, \left(1, v, \frac{|v|^2}{2}\right) \, dv \, = \, \bm{0}_{\RR^{5}}^\intercal.
			\]
			Plugging this expansion into \eqref{eqHydroGeneral}, we obtain the more detailed system
		  \begin{equation}
		    \label{eqHydroNonClosed}
				\left\{ \begin{aligned}
					& \partial_t \rho^\ve  + \diverg_x (\rho^\ve \, \bm{u}^\ve ) \,=\, 0, 
					\\
					& \,
					\\
					& \partial_t(\rho^\ve \, \bm{u}^\ve ) + \diverg_x \left (\rho^\ve \bm{u}^\ve \otimes \bm{u}^\ve + \rho^\ve T^\ve \left ( \bm{I} + \bm{\bar{A}}^\ve\right )\right ) \,=\, \bm{0}_{\RR^{3}}, 
					\\
					& \,
					\\
					& \partial_t E^\ve + \diverg_x \left (\frac 12 \rho^\ve |\bm{u}^\ve|^2 \bm{u}^\ve + \rho^\ve T^\ve \left (  \frac{3+2}{2}\bm{I} + \bm{\bar{A}}^\ve \right ) \bm{u}^\ve + \rho^\ve (T^\ve)^{3/2} \bm{\bar{B}}^\ve \right ) \,=\, 0,
				\end{aligned} \right.
			\end{equation}
			where the traceless matrix $\bm{\bar{A}}^\ve \in M_{3}$  and the vector $\bm{\bar{B}}^\ve\in \RR^{3}$ are given by
			\begin{equation}
			  \label{eqDefABV}
	      \left\{\begin{aligned}
		      & \bm{\bar{A}}^\ve := \frac{1}{\rho^\ve} \int_{\RR^{3}} \bm{A}(\bm{V}) f^\ve(v) \, dv, && \bm{A}(\bm{V}) = \bm{V} \otimes \bm{V} - \frac{|\bm{V}|^2}{3} \bm{I}, \\
		      & \bm{\bar{B}}^\ve := \frac{1}{\rho^\ve}\int_{\RR^{3}} \bm{B}(\bm{V}) f^\ve(v) \, dv, && \bm{B}(\bm{V}) = \frac12 \left [ |\bm{V}|^2 - (3+2)\right ] \bm{V}, 
	      \end{aligned}\right.
	    \end{equation}
	    and where we used the shorthand 
	    \[ \bm{V}(v) = \frac{v-\bm{u}}{\sqrt{T}}. \]
	    Depending on the order in $\ve$ of the truncation of the series \eqref{devChapEnsk}, we will obtain different hydrodynamic description of the fluid.
	    
			\medskip
			
			\subsubsection{Zeroth order: compressible Euler system}
			
			  At zeroth order with respect to  $\ve$, we have $f^\ve = \mathcal{M}_{\rho, \bm{u}, T}$. 
			  This distribution is in particular isotropic in $v-\bm{u}$ and its odd moments with respect to $(v-\bm{u})$ are all equal to zero. Since the matrix $\bm{\bar{A}}^\ve$ is traceless, we then have that 
			  \[ 
			    \bm{\bar{A}}_{Euler} := \frac{1}{\rho} \int_{\RR^{3}} \bm{A}(\bm{V}) \mathcal M_{\rho,\bm{u},T}(v) \,dv = \bm{0}_{M_{3}}.
			  \]
			  Moreover, since $ \bm{\bar{B}}^\ve$ involves odd, centered moments of $f^\ve$, we also obtain
			  \[
			    \bm{\bar{B}}_{Euler} := \frac{1}{\rho} \int_{\RR^{3}} \bm{B}(\bm{V}) \mathcal M_{\rho,\bm{u},T}(v) \,dv  = \bm{0}_{\RR^{3}}.
			  \]
			  Hence, the moments $(\rho, \bm{u}, T)$ are solution to the compressible Euler system \eqref{eqHydroClosedEuler}.
			  We notice in particular that the Maxwellian distribution in \eqref{devChapEnsk} is independent of $\ve$.
			  
			\medskip

			\subsubsection{First order: compressible Navier-Stokes system}
			
			  Going to the next order in $\ve$, we plug-in the expansion \eqref{devChapEnsk} in the Boltzmann equation \eqref{eqCollision}. 
			  Since the Maxwellian distribution is an equilibrium of the collision operator (according to \ref{hypEquilib}), the fluctuation $g^{(1)}$ is given by 
			  \begin{equation}
			    \label{eqCollisionExp}
			    \partial_t \mathcal M_{\rho,\bm{u},T} + v \cdot \nabla_x \mathcal M_{\rho,\bm{u},T} \,=\, \LL_{\mathcal M_{\rho,\bm{u},T}} \, g^{(1)}  + \mathcal{O}(\ve),
			  \end{equation}
			  where $(\rho, \bm{u}, T)$ are solution to the compressible Euler system \eqref{eqHydroClosedEuler}and $\LL_{\mathcal M}$ is the linearized\footnote{Namely the Frechet derivative of the collision operator.} collision operator around the Maxwellian distribution.
					  
        Besides, we also have
			  \begin{multline*}
			    \partial_t \mathcal M_{\rho,\bm{u},T} + v \cdot \nabla_x \mathcal M_{\rho,\bm{u},T} \, = \, \\
			      \mathcal M_{\rho,\bm{u},T} \left [ \partial_t \rho + v \cdot \nabla_x \rho + \frac{1}{\sqrt{T}} \left ( \bm{V} \cdot \partial_t \bm{u} + \bm{V}\otimes v : \nabla_x \bm{u}\right ) + \frac{1}{2 T} \left (|\bm{V}|^2 - 3\right ) \left (\partial_t T + v\cdot \nabla_x T\right ) \right ].
			  \end{multline*}
			  Then, using the conservation laws \eqref{eqHydroClosedEuler} in this equation to replace the time derivatives by only spatial ones, and dropping the terms of order $\ve$ in \eqref{eqCollisionExp}, we find after some computations that
			  \begin{equation}
			    \label{eqLinearG1}
			    \LL_{\mathcal M_{\rho,\bm{u},T}} g^{(1)} =  \mathcal M_{\rho,\bm{u},T}\left [ \bm{A}(\bm{V}):\bm{D}(\bm{u}) + 2 \bm{B}(\bm{V}) \cdot \nabla_x \sqrt{T}\right ],
			  \end{equation}
			  where $\bm{A}$, $\bm{B}$ and $\bm{V}$ are defined in \eqref{eqDefABV} and the traceless \emph{deformation tensor} $\bm{D}$ of $\bm{u}$ is given by		  
			  \[\bm{D}(\bm{u}) = \nabla_x \bm{u} + \left (\nabla_x \bm{u}\right )^\intercal - \frac{2}{3} \left (\diverg_x \bm{u}\right ) \bm{I}. \]
			  
			  Moreover, using the hypothesis \ref{hypConservations} on the conservation laws of the collision operator, it is possible to show that linear combinations of collisional invariants form exactly the kernel of the linear operator $\LL_{\mathcal M_{\rho,\bm{u},T}}$. In particular, we have the orthonormal family:
			  \[ 
			    \ker \LL_{\mathcal M_{\rho,\bm{u},T}} = \vect \left \{ \frac{1}{\rho},  \frac {\bm{V}}{\rho}, \frac{1}{2\rho} \left (|\bm{V}|^2 - 3 \right ) \right \}.
			  \]
			  Using the orthogonality properties of the moments of a Maxwellian distribution, we have that on $L^2 \left (\mathcal M_{\rho,\bm{u},T} \right )$,
			  \[ \bm{A} (\bm{V}), \bm{B}(\bm{V}) \perp \ker \LL_{\mathcal M_{\rho,\bm{u},T}}. \]
			  Since the operator $\LL_{\mathcal M_{\rho,\bm{u},T}}$ is  invertible on the orthogonal of its kernel and finally  using \eqref{eqLinearG1}, it yields
			  \begin{equation}
			    \label{eqG1}
			    g^{(1)} =  \LL_{\mathcal M_{\rho,\bm{u},T}}^{-1}\left ( \mathcal M_{\rho,\bm{u},T} \, \bm{A} \right ): \bm{D}(\bm{u}) \,+\, 2 \,\LL_{\mathcal M_{\rho,\bm{u},T}}^{-1}\left ( \mathcal M_{\rho,\bm{u},T} \, \bm{B} \right ) \cdot \nabla_x \sqrt{T}. 
			  \end{equation}
			  We can then plug this expression into the definition \eqref{eqDefABV} to obtain using some classical symmetry properties of the collision operator \cite{LevermoreMorokoffNadiga:98} that
			  \begin{equation}
			    \label{eqABNS}
			    \left\{\begin{aligned}
			      & \bm{\bar{A}}^\ve_{NS} := \frac{1}{\rho} \int_{\RR^{3}} \bm{A}(\bm{V}) \mathcal M_{\rho,\bm{u},T}(v)\left [1 + \ve \, g^{(1)}(v) \right ] dv = - \ve \frac{\mu}{\rho \, T} \bm{D}(\bm{u}), \\ 
	          & \bm{\bar{B}}^\ve_{NS} := \frac{1}{\rho}\int_{\RR^{3}} \bm{B}(\bm{V}) \mathcal M_{\rho,\bm{u},T}(v)\left [1 + \ve \, g^{(1)}(v) \right ] dv = - \ve \frac{\kappa}{\rho \,T^{3/2}} \nabla_x T.
	        \end{aligned}\right.
			  \end{equation}
	      The scalar quantities $\mu$ and $\kappa$ in \eqref{eqABNS}, respectively the \emph{viscosity} and the \emph{thermal conductivity}, are given by
	      \begin{gather*}
	        \mu := - T \int_{\RR^{3}} \mathcal M_{\rho,\bm{u},T}(v) \bm{A}(\bm{V}):\LL_{\mathcal M_{\rho,\bm{u},T}}^{-1}\left ( \mathcal M_{\rho,\bm{u},T} \, \bm{A} \right )(v) \, dv, \\
	        \kappa := - T \int_{\RR^{3}} \mathcal M_{\rho,\bm{u},T}(v) \bm{B}(\bm{V}) \cdot \LL_{\mathcal M_{\rho,\bm{u},T}}^{-1}\left ( \mathcal M_{\rho,\bm{u},T} \, \bm{B} \right )(v) \, dv.
	      \end{gather*}	      
	      They depend on the collision kernel of the model. For example, for the Boltzmann operator in the hard sphere case, it can be shown \cite{Golse:2005} that there exists some positive constants $\mu_0$, $\kappa_0$ such that
	      \[ \mu = \mu_0 \sqrt T \quad \text{and} \quad \kappa = \kappa_0 \sqrt{T}.\]
	      In the ES-BGK case, we have \cite{Struchtrup:2005}
	      \[ \mu = \frac{1}{1-\beta} \frac{\rho \, T}{\nu} \quad \text{and} \quad \kappa = \frac{5}{2} \frac{\rho \, T}{\nu}.\]
	      
		    Finally, the evolution of the macroscopic quantities at first order with respect to $\ve$ is given by the compressible Navier-Stokes equations
			  \begin{equation}
			    \label{eqHydroClosedNS}
					\left\{ \begin{aligned}
						& \partial_t \rho + \diverg_x (\rho\, \bm{u} ) \,=\, 0, 
						\\
						& \,
						\\
						& \partial_t(\rho \, \bm{u} ) + \diverg_x \left (\rho \, \bm{u} \otimes \bm{u} +  \rho \, T \, \bm{I} \right ) \,=\, \ve \diverg_x \left ({\mu} \,\bm{D}(\bm{u})\right ), 
						\\
						& \,
						\\
						& \partial_t E + \diverg_x \left ( \bm{u} \left ( E +\rho \, T\right ) \right ) \,=\, \ve \diverg_x \left ( \mu  \, \bm{D}(\bm{u}) \cdot \bm{u} + \kappa  \nabla_x T \right ).
					\end{aligned} \right.
				\end{equation}
				
				\begin{remark}
				  The matrix $\bm{\sigma} := -\mu  \, \bm{D}(\bm{u})$ is sometimes called \emph{viscosity tensor} and the vector $\bm{q} := -\kappa  \nabla_x T$ is the \emph{heat flux}.
				\end{remark}

			\subsubsection{Second order: Burnett equations}
				
				Pushing the expansion \eqref{devChapEnsk} at second order in $\ve$, we can use the same type of argument that for the compressible Navier-Stokes system to obtain another correction of the compressible Euler equations: the Burnet system.
				Although this system is ill-posed \cite{Golse:2005}, the computation of its coefficients is still possible.
				We have for the BGK case $\beta = 0$  \cite{Struchtrup:2005}:
				\begin{align}
				  \bm{\bar{A}}^\ve_{Burnett} & := \frac{1}{\rho}\int_{\RR^{3}} \bm{A}(\bm{V}) \mathcal M_{\rho,\bm{u},T}(v)\left [1 + \ve \, g^{(1)}(v) + \ve^2  g^{(2)}(v) \right ] dv \notag \\
				    & = - \ve \frac{\mu}{\rho \,T} \bm D( \bm{u} ) 
				        - 2 \ve^2 \frac{\mu^2}{\rho^2 T^{2}}\bigg \{ -\frac{T}{\rho} {\rm Hess}_x(\rho ) + \frac{T}{\rho^2} \nabla_x \rho \otimes \nabla_x \rho - \frac{1}{\rho} \nabla_x T \otimes \nabla_x \rho \notag \\
				    &  \qquad \qquad \qquad \qquad \qquad \qquad \ + \left (\nabla_x \bm{u} \right )  \left (\nabla_x \bm{u} \right )^\intercal - \frac{1}{3}\bm D( \bm{u} ) \diverg_x( \bm{u} ) + \frac{1}{T} \nabla_x T \otimes \nabla_x T\bigg \};  \label{eqAepsBurnett} \\				  
          \bm{\bar{B}}^\ve_{Burnett} & := \frac{1}{\rho}\int_{\RR^{3}} \bm{B}(\bm{V}) \mathcal M_{\rho,\bm{u},T}(v)\left [1 + \ve \, g^{(1)}(v) + \ve^2  g^{(2)}(v) \right ] dv \notag \\
				    & = - \ve \frac{\kappa}{\rho \,T^{3/2}} \nabla_x T 
				        - \ve^2 \frac{\mu^2}{\rho^2 T^{5/2}}\bigg \{ +\frac{25}{6} \left (\diverg_x \bm{u}\right ) \nabla_x T \notag \\
				    &  \qquad \qquad \qquad \qquad \qquad \qquad \quad \ \, - \frac{5}{3} \left [T \diverg_x \left (\nabla_x \bm{u}\right ) + \left (\diverg_x\bm{u} \right ) \nabla_x T + 6 \left (\nabla_x \bm{u} \right )  \nabla_x T\right ] \notag \\
				    &  \qquad \qquad \qquad \qquad \qquad \qquad \quad \ \, + \frac{2}{\rho} \bm D( \bm{u} ) \, \nabla_x \left (\rho \, T \right ) + 2 \, T \diverg_x \left (\bm D( \bm{u} )\right ) + 16 \bm D( \bm{u} ) \,\nabla_x T \bigg \} .
            \label{eqBepsBurnett}
				\end{align}

	 	\subsection{From Fluid to Kinetic: the Moment Realizability Criterion}
	 	  \label{subFluid2Kin}

			The matrix $\bm{\bar{A}}^\ve$ and the vector $\bm{\bar{B}}^\ve$ will allow us to define our hydrodynamic break down criterion.  Let us set the vector of the reduced collisional invariants  for $\bm{V}=(v-\bm{u^\ve})/\sqrt{T^\ve}$,
      \[
		    \bm{m} := \left (1, \bm{V}, \left (\frac{2}{3}\right )^{1/2}\left (\frac{|\bm{V}|^2}{2} - \frac{3}{2}\right ) \right ).
		  \]  
		  We then define the so-called \emph{moment realizability matrix} by setting
		  \begin{equation}
		    \label{eqMomRealiz}
		    \bm{M} := \frac{1}{\rho^\ve}\int_{\RR^{3}} \bm{m} \otimes \bm{m} \, f^\ve(v) \, dv.
		  \end{equation}
		  By using the orthogonality properties of the moments of a Maxwellian distribution and \eqref{devChapEnsk}, we have
	    \begin{align*}
	     \bm{M} & = \frac{1}{\rho^\ve} \int_{\RR^d} 
	     \begin{pmatrix}
	       1 & \bm{V}^\intercal & \left (\frac{2}{3}\right )^{1/2}\left (\frac{|\bm{V}|^2}{2} - \frac{3}{2}\right ) \\
	       \bm{V} & \bm{V} \otimes \bm{V} & \left (\frac{2}{3}\right )^{1/2}\left (\frac{|\bm{V}|^2}{2} - \frac{3}{2}\right ) \bm{V} \\
	       \left (\frac{2}{3}\right )^{1/2}\left (\frac{|\bm{V}|^2}{2} - \frac{3}{2}\right ) & \left (\frac{2}{3}\right )^{1/2}\left (\frac{|\bm{V}|^2}{2} - \frac{3}{2}\right ) \bm{V}^\intercal & \frac{2}{3} \left (\frac{|\bm{V}|^2}{2} - \frac{3}{2}\right )^2 
	     \end{pmatrix} f^\ve(v) \, dv \notag\\
	     & =
	     \begin{pmatrix}
	       1                  & \bm{0}_{\RR^{3}}^\intercal                 & 0 \\
	       \bm{0}_{\RR^{3}} & \bm{I} +\bm{\bar{A}}^\ve                     & \left (\frac{2}{3}\right )^{1/2} \bm{\bar{B}}^\ve \\
	       0 & \left (\frac{2}{3}\right )^{1/2} (\bm{\bar{B}}^\ve)^\intercal & \bar{C}^\ve
	     \end{pmatrix},   
	    \end{align*}
		  where $\bar{C}^\ve$ is the dimensionless fourth order moment of $f^\ve$:
		  \[
		     \bar{C}^\ve := \frac{2}{3 \rho}\int_{\RR^{3}} \left [\frac{|\bm{V}|^2}{2} - \frac{3}{2}\right ]^2 f^\ve(v) \, dv.
		  \]
		  For the sake of simplicity, let us introduce the change of basis $\bm Q$ by setting
		  \[
		    \bm{Q} := 
		    \begin{pmatrix}
	       1                  & \bm{0}_{\RR^{3}}^\intercal                  & 0 \\
	       \bm{0}_{\RR^{3}} & \bm{I}                                        & \bm{0}_{\RR^{3}} \\
	       0 &\displaystyle  - \left (\frac{2}{3}\right )^{2/3} \frac{(\bm{\bar{B}}^\ve)^\intercal}{\bar{C}^\ve} & 1
	     \end{pmatrix}.
		  \]
		  Then we have the following relation for $\bm{M}$:
		  \[ \bm M =
		    \bm{Q}^{-1} \begin{pmatrix}
		      1 & \bm{0}_{\RR^{3}}^\intercal       & 0 \\
		       \bm{0}_{\RR^{3}}  & \displaystyle\bm{I} + \bm{\bar{A}}^\ve - \frac{2}{3\,\bar C^\ve} \bm{\bar{B}}^\ve \otimes \bm{\bar{B}}^\ve &  \bm{0}_{\RR^{3}}  \\
		      0 &  \bm{0}_{\RR^{3}}^\intercal                & \bar C^\ve
		    \end{pmatrix} \bm{Q}.
		  \]
		  Since $\bar C^\ve$ is a nonnegative quantity and $\bm{M}$ is by construction a positive definite matrix, the matrix
  		\begin{equation}
		    \label{eqMomRealizReduced}
		    \mathcal{V} := \bm{I} + \bm{\bar{A}}^\ve - \frac{2}{3\,\bar C^\ve} \bm{\bar{B}}^\ve \otimes \bm{\bar{B}}^\ve
		  \end{equation}
		  is also a positive definite matrix. 

      Now from these remarks, let us define a criterion to determine the appropriate model -- fluid or kinetic -- to be used.

On the one hand, consider the zeroth order model with respect to $\ve$, that is the compressible Euler system. If we truncate the expansion at first order in $\ve$, we get that $\bm{\bar{A}}_{Euler} = \bm{0}_{M_{3}}$ and $\bm{\bar{B}}_{Euler} = \bm{0}_{R^{3}}$. Moreover, we also have in that case  $\bar{C}^\ve = 1$ and then
				\[	
					\mathcal{V}_{Euler} := \mathcal{V}_{1} = \bm{I}. 
				\]

 On the other hand, consider  the first order model, that is, the  compressible Navier-Stokes system. By cutting the Chapman-Enskog expansion \eqref{devChapEnsk} at the first order with respect to $\ve$ (\ie Navier-Stokes order), we can compute explicitly the matrix $\mathcal{V}_{\ve}$. 
		    We have in this case using the expressions \eqref{eqABNS} and by symmetry arguments \cite{LevermoreMorokoffNadiga:98}, that $\bar C^\ve = 1$, hence
		    \begin{equation}
		      \label{eqMomRealizNS}
		      \mathcal{V}_{NS} := \mathcal{V}_{\ve} = \bm{I} - \ve \frac{\mu}{\rho T} \bm{D}\left (\bm{u}\right ) - \ve^2 \frac{2}{3} \frac{\kappa^2}{\rho^2 T^3} \nabla_x T \otimes \nabla_x T,
		    \end{equation}
		    where $(\rho, \bm{u}, T)$ are solution to the Navier-Stokes equations \eqref{eqHydroClosedNS}.

 Hence, we claim that the compressible Euler system is correct when the matrix $ \mathcal{V}_{NS}$ behaves like the matrix $\mathcal{V}_{Euler}={\bm I}$, that is, it is  positive definite and if \emph{its eigenvalues are close to $1$ or not}: The Euler description of the fluid will be considered incorrect if 
				\begin{equation}
			  	\label{crit:CompEuler}
			  	\left |\lambda_{NS}-1\right | > \eta_0, \qquad\forall \lambda_{NS} \in\, {\rm Sp}(\mathcal{V}_{NS}),
				\end{equation}
	where $\eta_0$ is a small parameter (here we take $\eta_0=10^{-2}$).

	 More generally,  we denote by $f^\ve_k$ the $k^\text{th}$ order truncation of the Chapman-Enskog expansion \eqref{devChapEnsk}:
      \begin{equation}
        \label{devChapEnsk_k} 
        f^\ve_k := \mathcal M_{\rho,\bm{u},T}\left [1 + \ve \, g^{(1)} + \ve^2 \, g^{(2)} + \ldots + \ve^k  \, g^{(k)} \right ].
      \end{equation}

For a given truncation (\ref{devChapEnsk_k}) of order $k$, we will say that the fluid model associated is \emph{incorrect}  at point $(t,x)$ if we have
	        \begin{equation}
	          \label{crit:orderK}
		  		  \left |\lambda_{\ve^k}-\lambda_{\ve^{k+1}}\right | > \eta_0, \quad\forall \, \lambda_{\ve^k} \in\, {\rm Sp}(\mathcal{V}_{\ve^k}), \ \lambda_{\ve^{k+1}} \in\, {\rm Sp}(\mathcal{V}_{\ve^{k+1}}).
				  \end{equation}

		\subsection{From Kinetic to Fluid}
		  \label{subKin2Fluid}
	    
	    Knowing the full kinetic description of a gas, there exists a large number of methods \cite{Saint-Raymond:2009} to decide how far this gas is from the thermal equilibrium, \ie the fluid regime. 
	    We decide to use a simple comparison between the kinetic density $f^\ve$, solution to the collisional equation \eqref{eqCollision} and the truncated Chapman-Enskog distribution $f^\ve_k$ given by \eqref{devChapEnsk_k}, whose moments match the one of $f^\ve$, and whose order $k$ corresponds to the order of the macroscopic model considered.
	    
%
%
%
	    

     Our criterion is then the following: The kinetic description at point $(t,x)$ corresponds to an \emph{hydrodynamic closure of  order $k$} if
		  \begin{equation}
		    \label{eqKintoFluid1}
		    \left \| f^\ve(t,x,\cdot) - f^\ve_k(t,x,\cdot) \right \|_{L^1_v} \leq \delta_0,
		  \end{equation}
		  where $\delta_0$ is a small parameter (we take $\delta_0 = 10^{-4}$).

		  \begin{example}
		    If $k=1$ (Compressible Euler setting), this criterion corresponds to the natural one
	      \[
	        \|f^\ve(t,x,\cdot) - \mathcal M_{\rho,\bm{u},T}^\ve(t,x,\cdot)\|_{L^1} \leq \delta_0,
	      \]
	      namely to check if the system is locally at the thermodynamic equilibrium or not.
	    \end{example}

		  \begin{remark}		  
		    In particular, if we perform the Chapman-Enskog expansion \eqref{devChapEnsk} of $f^\ve$, the criterion $\eqref{eqKintoFluid1}$ corresponds to the fact that the remainder term in this expansion is small in $L^1$ norm, because it is then given by
		    \[
		      \left \| \sum_{n > k} \ve^i \, g^{(i)}(t,x,\cdot)\right \|_{L^1_v} \leq \delta_0.
		    \]
		  \end{remark}
	    
Moreover,  for numerical purposes, it can be interesting to take an additional criterion on  
\begin{equation}
\label{eqKintoFluid2}
\frac{\Delta t}{\ve} \gg 1,
  \end{equation}  
where $\Delta t$ is the time step.
	      Indeed, the relaxation time of equation \eqref{eqCollision} toward the Maxwellian distribution is of order $\ve/\Delta t$. Hence for small $\ve$ or large time step $\Delta t$ the solution is at thermodynamical equilibrium.

	\section{Numerical schemes}
	  \label{secNumSim}    
			  
		  \subsection{Systems of Conservation Laws}

			  In this subsection, we shall focus on the space discretization of the system of $n$ conservation laws
			  \begin{equation}
				  \label{sysConservLaws}
				  \left \{ \begin{aligned}
				   & \frac{\partial u}{\partial t} + \diverg_x F(U) = 0, \ \forall \, (t,x) \in \RR_+ \times \Omega, \\
				   & \, \\
				   & u(0,x) = u_0(x), 
				   \end{aligned} \right.
				\end{equation}
				for a smooth function $F : \RR^n \to M_{n\times d_x}(\RR)$ and a Lipschitz-continuous domain $\Omega \subset \RR^{d_x}$. 
			
				Here we apply  finite volume schemes using central Lax Friedrichs schemes with slope limiters (see \eg Nessyahu and Tadmor \cite{NessyahuTadmor:1990}).  
				
			\medskip

\subsection{ES-BGK Equation}
		\label{subNumSchemes}		
			  We now focus briefly on the time evolution of the ES-BGK equation
			  \begin{equation} 
					\label{eqESBGK}
					\left\{ \begin{aligned}
					  & \frac{\partial f^\ve}{\partial t} + v \cdot \nabla_x f^\ve \,=\, \frac{\nu}{\ve} \left( \G[f] \,-\, f \right), 
					  \\
				  	&  \;
					  \\
					  & f^\ve(0, x, v) = f_{0}(x,v),
					\end{aligned} \right.
				\end{equation}
				We adopt the approach of Filbet and Jin \cite{FilbetJin:2010}, that is we discretize the time using a first order Implicit-Explicit (\emph{IMEX}) scheme.
				Since the convection term in \eqref{eqESBGK} is not stiff, we treat it explicitly, and we use an implicit solver only for the stiff source term on the right hand side.

    \subsection{Evolving the Variables and Coupling the Equations}
      \label{subEvoCouple}  
      
		  We are now interested in evolving in time the hybrid scheme. Let us consider the case of a fluid closure of order $k$.
		  At a given time $t^n$, we denote by $K_i$ a control volume, the space domain $\Omega = \Omega_f^n \sqcup \Omega_k^n$ is decomposed in
		  \begin{itemize}
		    \item \emph{Fluid cells} $K_i \subset \Omega_f^n$, described by the hydrodynamic fields
		      \[
		        U_i^n := \left (\rho_i^n, \bm{u}_i^n, T_i^n\right ) \simeq \left (\rho(t^n,x_i), \bm{u}(t^n,x_i), T(t^n,x_i)\right ) ; 
		      \]
		    \item \emph{Kinetic cells} $K_j \subset \Omega_k^n$, described by the particle distribution function
		    \[
		      f_j^n(v) \simeq f(t^n,x_j,v), \quad \forall v \in \RR^{3}.
		    \]
		  \end{itemize}		  
		  The evolution of the whole system depends on the type of cell we consider. The algorithm used is the following:
		  
		  \begin{itemize}
		    \item In a fluid cell $K_i \subset \Omega_f^n$, compute the eigenvalues of the reduced  moment realizability matrix $\mathcal{V}_{\ve^k}$:
		    
		      \begin{itemize}
		        \item If the criterion \eqref{crit:orderK} is wrong, evolve the fluid equations at point $x_i$ with initial condition $U_i^n$  to obtain $U_i^{n+1}$;
		        \smallskip
		        \item In the other case, the regime is no longer fluid but kinetic, then ``lift'' the macroscopic fields into the kinetic grid, by taking for new distribution $f_i^n$ a Maxwellian\footnote{The proper way to do so since the velocity space is discrete is to consider discrete velocity Maxwellians, as introduced by Berthelin, Tzavaras and Vasseur in \cite{BerthelinTzavarasVasseur:2009}.}, whose moments are given by $U_i^{n}$:
		          \[
		            f_i^n(v) := \mathcal{M}_{\rho_{i}^n, \bm{u}_{i}^n, T_{i}^n}(v), \quad \forall v \in \RR^{3}.
		          \]
		          Evolve the kinetic equation at point $x_i$ with initial condition $f_i^n$ to obtain $f_i^{n+1}$;
		        \smallskip  
		        \item Set $\Omega_f^{n+1} := \Omega_f^n\setminus K_i$ and $\Omega_k^{n+1} := \Omega_k^n \cup K_i$.
		      \end{itemize}
		      
		    \item In a kinetic cell $K_j \subset \Omega_k^n$, evaluate the criteria \eqref{eqKintoFluid1}-\eqref{eqKintoFluid2}:
		    
		      \begin{itemize}
		        \item If both are correct, evolve the kinetic equation at point $x_j$ with initial condition $f_i^n$ to obtain $f_j^{n+1}$;
		        \smallskip
		        \item In the other case, the regime is fluid, then project the kinetic distribution towards the macroscopic fields, by setting
		          \[
		            U_{j}^n := \int_{\RR^d} f_{j}^n \, \varphi(v) \, dv, \quad \varphi(v) = \left (1, v, \frac{1}{3 \rho^n_j }|v - \bm{u_j^n}|^2 \right ).
		          \]  
		          Evolve the fluid equation at point $x_j$  with initial condition $U_j^n$  to obtain $U_j^{n+1}$;
		        \smallskip
		        \item Set $\Omega_k^{n+1} := \Omega_k^n\setminus K_j$ and $\Omega_f^{n+1} := \Omega_f^n \cup K_i$.
		      \end{itemize}
		    
		  \end{itemize}
			  
\begin{figure}[h!]
        \psset{xunit=1.cm,yunit=1.cm,algebraic=true,dotstyle=o,dotsize=4pt 0,linewidth=0.8pt,arrowsize=3pt 2,arrowinset=0.25}
        \begin{center}
					\begin{pspicture*}(-3.,0.3)(11.,3.2)
						\psline(1,1)(4,1)
						\psline(4,1)(7,1)
						\psline[linestyle=dashed,dash=2pt 2pt](4,0.5)(4,2.7)					
						\psline(7,1)(10,1)
						\psline(1,1)(-2,1)
						\rput[tl](-0.8,0.8){$K_{i-2}$}
						\rput[tl](2.3,0.8){$K_{i-1}$}
						\rput[tl](5.3,0.8){$K_{i}$}
						\rput[tl](8.3,0.8){$K_{i+1}$}
						\rput[tl](5.3,2.5){$f_i^n$}
						\rput[tl](-0.8,2.5){$U_{i-2}^n$}
						\rput[tl](2.3,2.5){$U_{i-1}^n$}
						\rput[tl](8.3,2.5){$f_{i+1}^n$}
						\begin{scriptsize}
							\rput[tl](6.5,1.75){Kinetic}
							\rput[tl](0.,1.75){Hydrodynamic}
						\end{scriptsize}
						\psdots[dotstyle=+,linecolor=blue](1,1)
						\psdots[dotstyle=+,linecolor=blue](4,1)
						\psdots[dotstyle=+,linecolor=blue](7,1)
						\psdots[dotstyle=+,linecolor=blue](10,1)
						\psdots[dotstyle=+,linecolor=blue](-2,1)
					\end{pspicture*}
			  \end{center}
			  \caption{Transition between fluid and kinetic cells.}
			  \label{figInterfaceHydroKin}
			\end{figure}
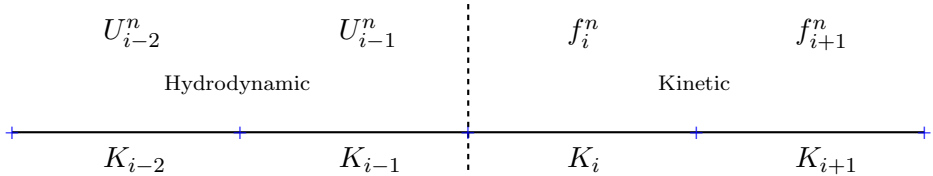

			It now remains to consider what happens between two cells of different types. Consider the situation at time $t^n$ where the cells $K_{i-2}$ and $K_{i-1}$ are fluid, and the cells $K_i$ and $K_{i+1}$ are kinetic (as described in Figure \ref{figInterfaceHydroKin}).
			  
			  \begin{itemize}
			    \item To evolve $f_i^n$ with a finite volume method as described in Section \ref{subNumSchemes}, a stencil of two ghost cells is needed on the left. Since we don't have boundary conditions prescribed by the problem between cells, we lift the  hydrodynamic fields by setting
			      \[
			        f_{i-1}^n(v) := \mathcal{M}_{\rho_{i-1}^n, \bm{u}_{i-1}^n, T_{i-1}^n}(v), \quad f_{i-2}^n(v) := \mathcal{M}_{\rho_{i-2}^n, \bm{u}_{i-2}^n, T_{i-2}^n}(v), \quad \forall v \in \RR^{3};
			      \]
			    			    
			    \item To evolve $U_{i-1}^n$ with a finite volume method, a {stencil} of two ghost cells is needed on the right. Similarly, we project the kinetic density by setting
			      \[
			        U_{i}^n := \int_{\RR^d} f_{i}^n \left (1, v, \frac{1}{3 \rho^n_i }|v - \bm{u}_i^n|^2\right ) dv, \quad U_{i+1}^n := \int_{\RR^d} f_{i+1}^n \left (1, v, \frac{1}{3 \rho^n_{i+1} }|v - \bm{u}_{i+1}^n|^2\right )  dv,
			      \]
			      and we match the hydrodynamic fluxes with the kinetic ones.

			  \end{itemize}
			  
			\begin{remark}[Limitations of this approach] 
			  Each zone must be at least as wide as the stencil, and the extensions to non-cartesian grids seems nontrivial.
			\end{remark}
			  
\section{Numerical Simulations}
	    \label{subNumSim}
	  
	    We take in all the simulations $d_x = 1$. 
	    In particular,  the moment realizability matrices $\bm{V}_{\ve^k}$ are diagonal. In the \eg Navier-Stokes case, it is  given by
	    \[ 
	      \mathcal{V}_{NS} = \begin{pmatrix}
	                      \displaystyle 1 - \ve \frac{\mu}{\rho T} \partial_x u^x - \ve^2\frac{\kappa^2}{\rho^2 T^3} (\partial_x T)^2  & 0 & 0\\
	                      0 & \displaystyle 1 + \ve \frac{\mu}{\rho T} \partial_x u^x & 0 \\
	                      0 & 0 & \displaystyle 1 + \ve \frac{\mu}{\rho T} \partial_x u^x \\
	                    \end{pmatrix},
	    \]
	    where $\bm{u}=(u^x,u^y,u^z)$.
	    We can then read its eigenvalues on its diagonal. 
	    The criterion for $k=0$ for a fluid cell to be kinetic at the next iteration is then
	    \begin{equation}
	      \label{criterion1dx3dv}
	      \left |\ve \frac{\mu}{\rho T} \partial_x u^x + \ve^2\frac{\kappa^2}{\rho^2 T^3} (\partial_x T)^2\right |  \leq \eta_0 
	      \quad \text{ or } \quad \left |\ve \frac{\mu}{\rho T} \partial_x u^x\right | \leq \eta_0.
	    \end{equation}
	    Using the expression of the Burnett coefficients \eqref{eqAepsBurnett}-\eqref{eqBepsBurnett}, we  can easily write the same type or criterion for the Navier-Stokes closure $k=1$.

\subsection{Test 1: Riemann problem}
	   
This test deals with the numerical solution of the non homogeneous $1D\times 3D$ BGK equation  \eqref{eqOpESBGK}. We present some results for one dimensional Riemann problem and compare them with the numerical solution obtained by solving the full kinetic equation on a fine mesh. We have computed an approximation for different Knudsen numbers from rarefied regime up to the fluid limit and report the results for $\ve=10^{-2}$ and $10^{-3}$. 

\begin{figure}
  \includegraphics{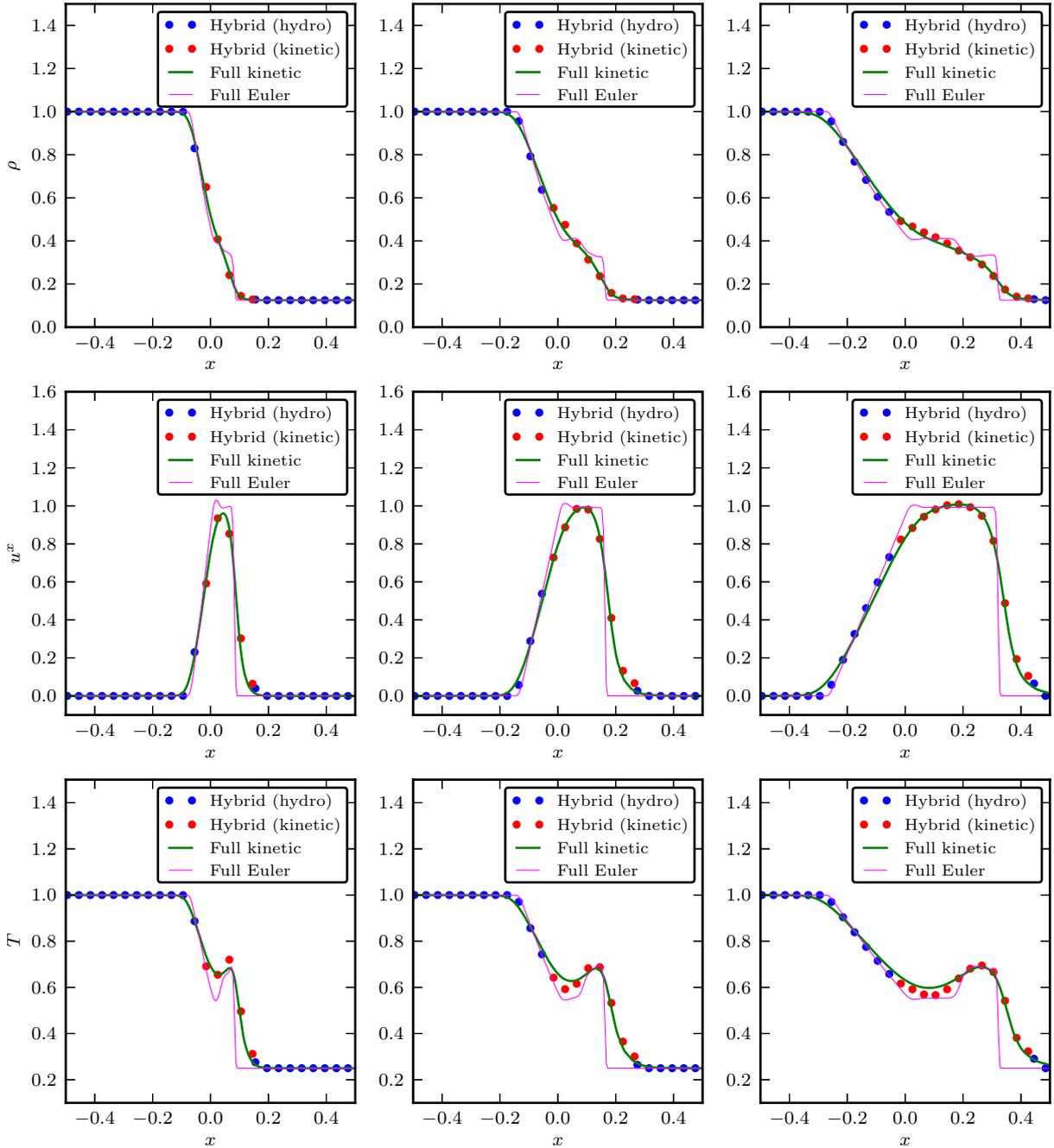}
  \caption{\textbf{Test 1 - Riemann problem with $\ve = 10^{-2}$ :} \emph{Order 0 (Euler)}; Density, mean velocity and temperature at times  $t = 0.05$, $0.10$ and $0.20$.}
  \label{figSodE2-Euler}
\end{figure}

\begin{figure}
  \includegraphics{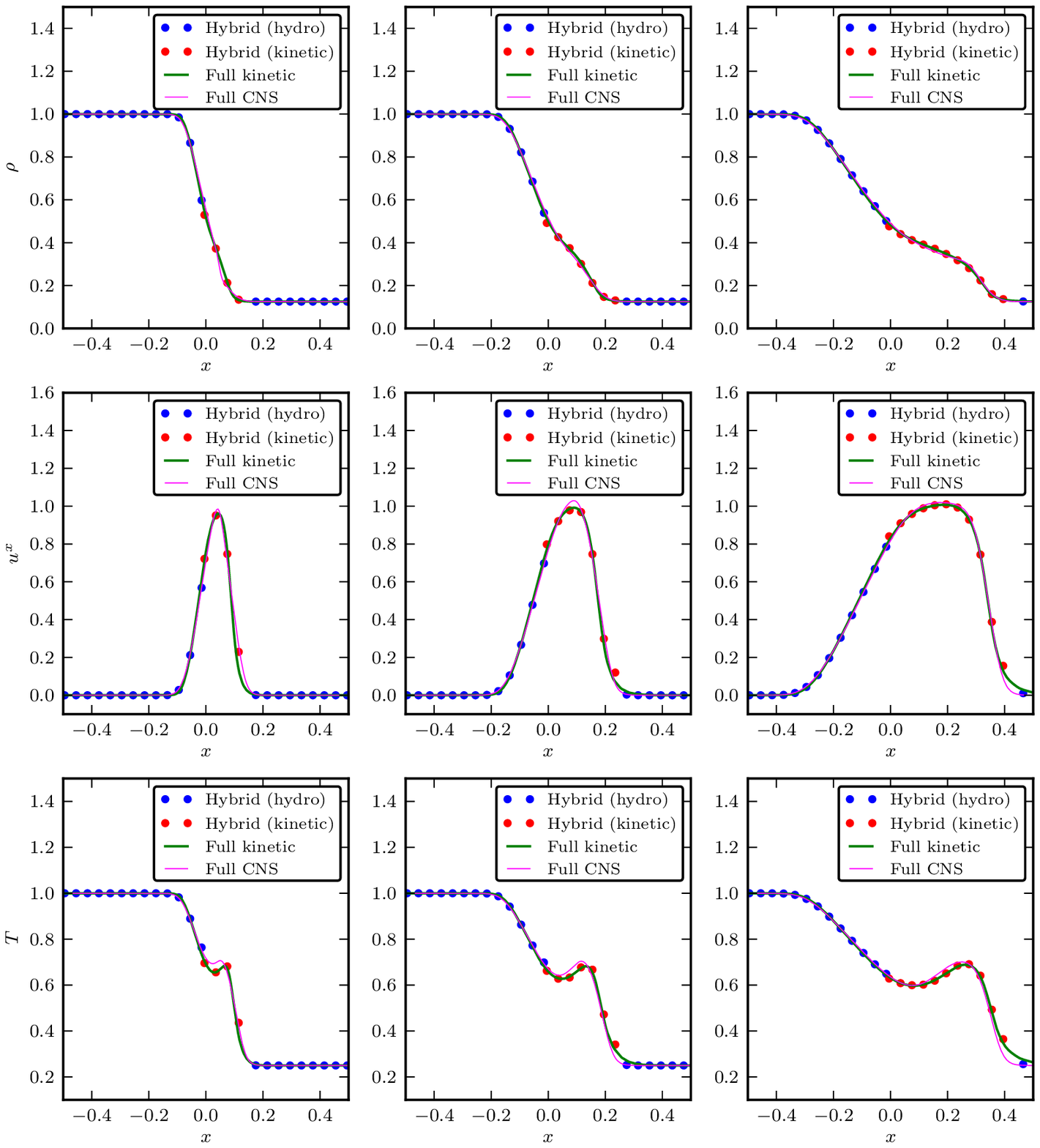}
	\caption{\textbf{Test 1 - Riemann problem with $\ve = 10^{-2}$ :} \emph{Order 1 (CNS)}; Density, mean velocity and temperature at times  $t = 0.05$, $0.10$ and $0.20$.}
	\label{figSodE2-CNS}
\end{figure}

\begin{figure}
  \includegraphics{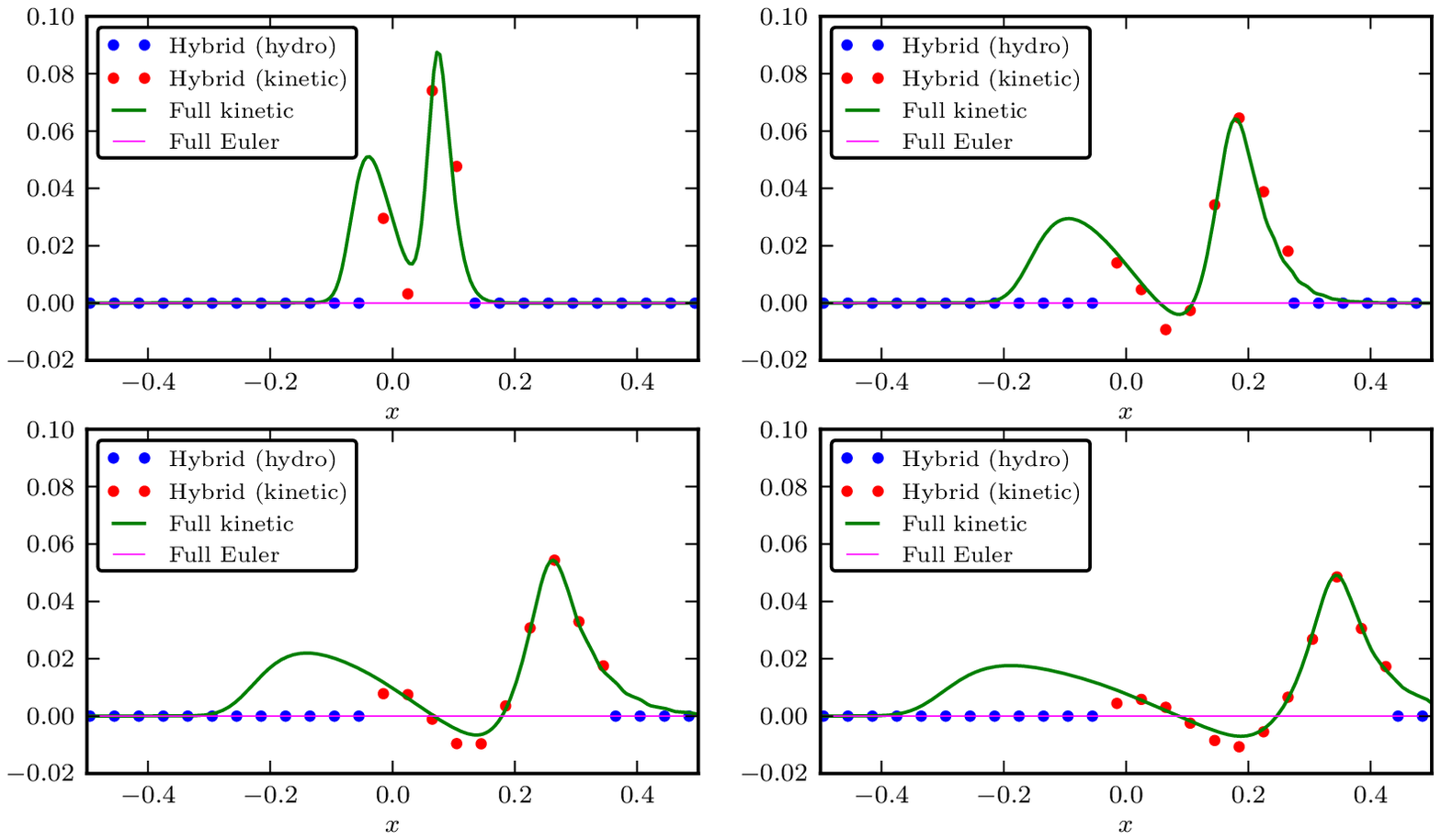}
  \caption{\textbf{Test 1 - Riemann problem with $\ve = 10^{-2}$ :} \emph{Order 0 (Euler)}; heat flux at times  $t = 0.05$, $0.10$, $0.15$ and $0.20$.}
  \label{figSodE2_q-Euler}
\end{figure}

\begin{figure}
  \includegraphics{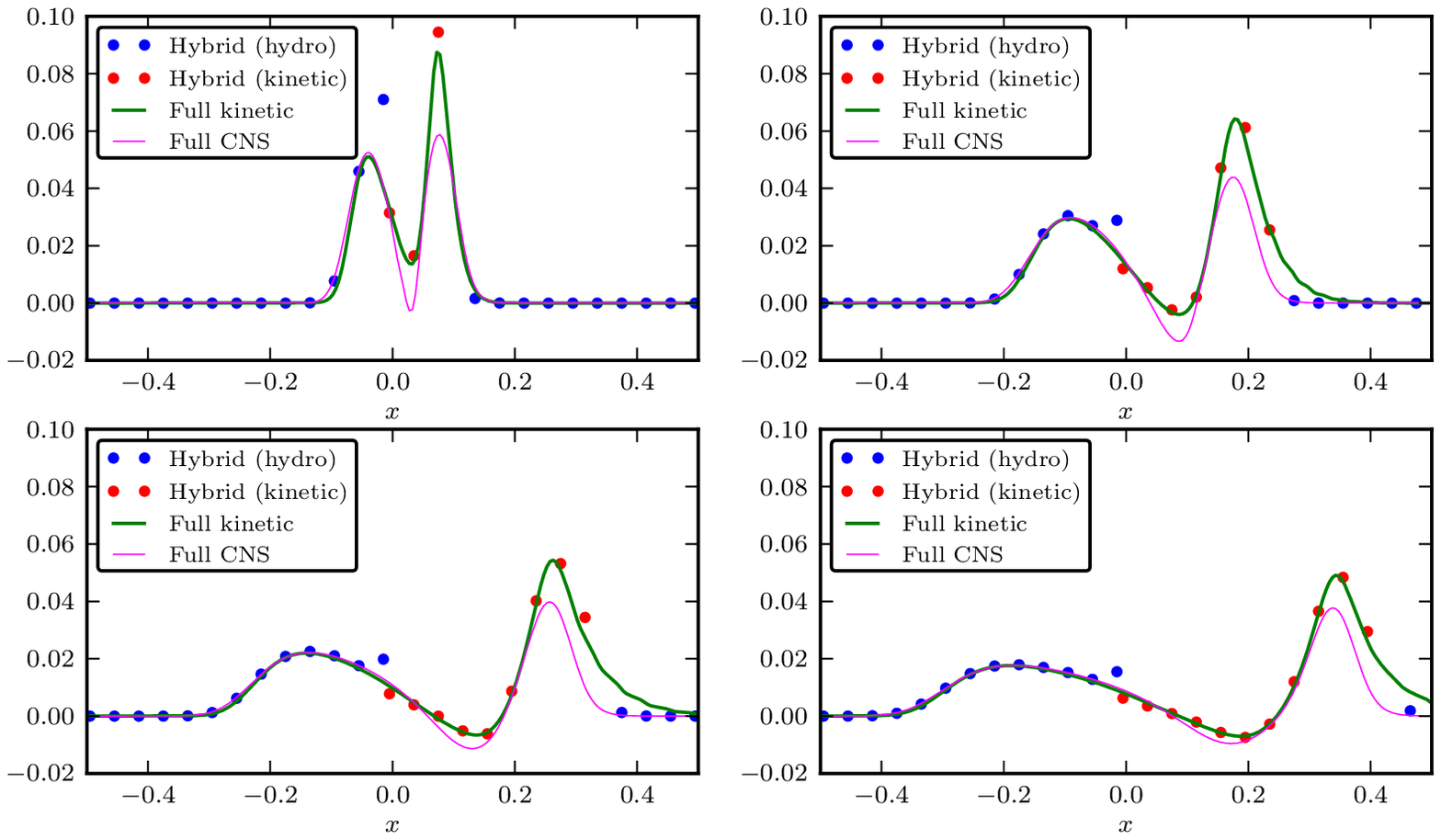}
	\caption{\textbf{Test 1 - Riemann problem with $\ve = 10^{-2}$ :} \emph{Order 1 (CNS)}; heat flux at times  $t = 0.05$, $0.10$, $0.015$ and $0.20$.}
	\label{figSodE2_q-CNS}
\end{figure}

More precisely, the initial data is given by
\[
 f^{in}(x,v) = \mathcal{M}_{\rho(x), \bm{u}(x), T(x)}(v), \quad \forall x\in[-0.5,0.5], \quad v\in [-8,8]^3,
 \]
with
\[
 \left (\rho(x), \bm{u}(x), T(x)\right ) = \left \{ 
\begin{aligned}
			            & (1,0,0,0,1) && \text{ if } x < 0, \\   
			            & (0.125,0,0,0,0.25) && \text{ if } x \geq 0.
\end{aligned}
 \right.;
\]
On the one hand, in Figures \ref{figSodE2-Euler} and \ref{figSodE2_q-Euler}, we plot the results obtained in the rarefied regime with $\ve=10^{-2}$, for the zeroth order model, namely the Euler dynamics.  
The kinetic reference solution is computed  with $200\times 128\times 32\times21$ cells in phase space, the fluid reference solution with $200$ points whereas the hybrid scheme is used with $100$ points in $x$ and the size of the velocity grid is $32 \times 32 \times 32$ points. 
We observe that the fluid solution is far from the kinetic one, which was expected since the Knudsen number is large. Nevertheless, the hybrid scheme behaves very nicely in this case, detecting correctly the non-equilibrium zone and the solution is close to the kinetic one. This error is mainly due to the application of the Euler equations for which the heat flux is zero (hence some errors in this particular quantity, see Figure \ref{figSodE2_q-Euler}).
Then, in Figures \ref{figSodE2-CNS} and \ref{figSodE2_q-CNS}, we perform the same simulations for the first order, Compressible Navier-Stokes (CNS) model.
Although the fluid solution is still far from the kinetic one, we observe that the result of the kinetic solver is in almost perfect agreement with the reference solution, even in large time. This can also be observed in the values of the heat flux, which are close to the reference ones. 

%
%
%
%

Concerning the computational times for the same configuration ($N_x=100$ and $N_v=32^3$), the hybrid schemes, both zeroth and first order, are more efficient than  the kinetic models since the computational time are respectively 1.9  and 4.4 times faster, even for such large values of $\ve$, corresponding to the rarefied case.
More details can be found in Table \ref{tabCompTime}.

\begin{table}

  \begin{tabular}{|c|c|c|c|c|c|}
    \hline
    Test & Sod $10^{-2}$ & Sod $10^{-3}$ & Blast $10^{-2}$ & Blast $5 \cdot 10^{-3}$ & Blast $10^{-3}$ \\
    \hline 
    Euler             & 0.03 & 0.03 & 0.02 & 0.02 & 0.02  \\
    \hline 
    CNS               & 0.08 & 0.09 & 0.1  & 0.1  & 0.11  \\
    \hline 
    BGK               & 113  & 120  & 160  & 161  & 158   \\
    \hline 
    Hybrid (Euler)    & 61.2 & 20.1 & 57   & 11   & 0.12  \\
    \hline 
    Hybrid (CNS)      & 25.6 & 4.9  & 23   & 18   & 3.3   \\
    \hline
  \end{tabular}
  \caption{Comparison of the computational times (sec), $t=0.10$, $N_x = 100$, $N_v = 32^3$.}
  \label{tabCompTime}
  
\end{table}

\begin{figure}
  \includegraphics{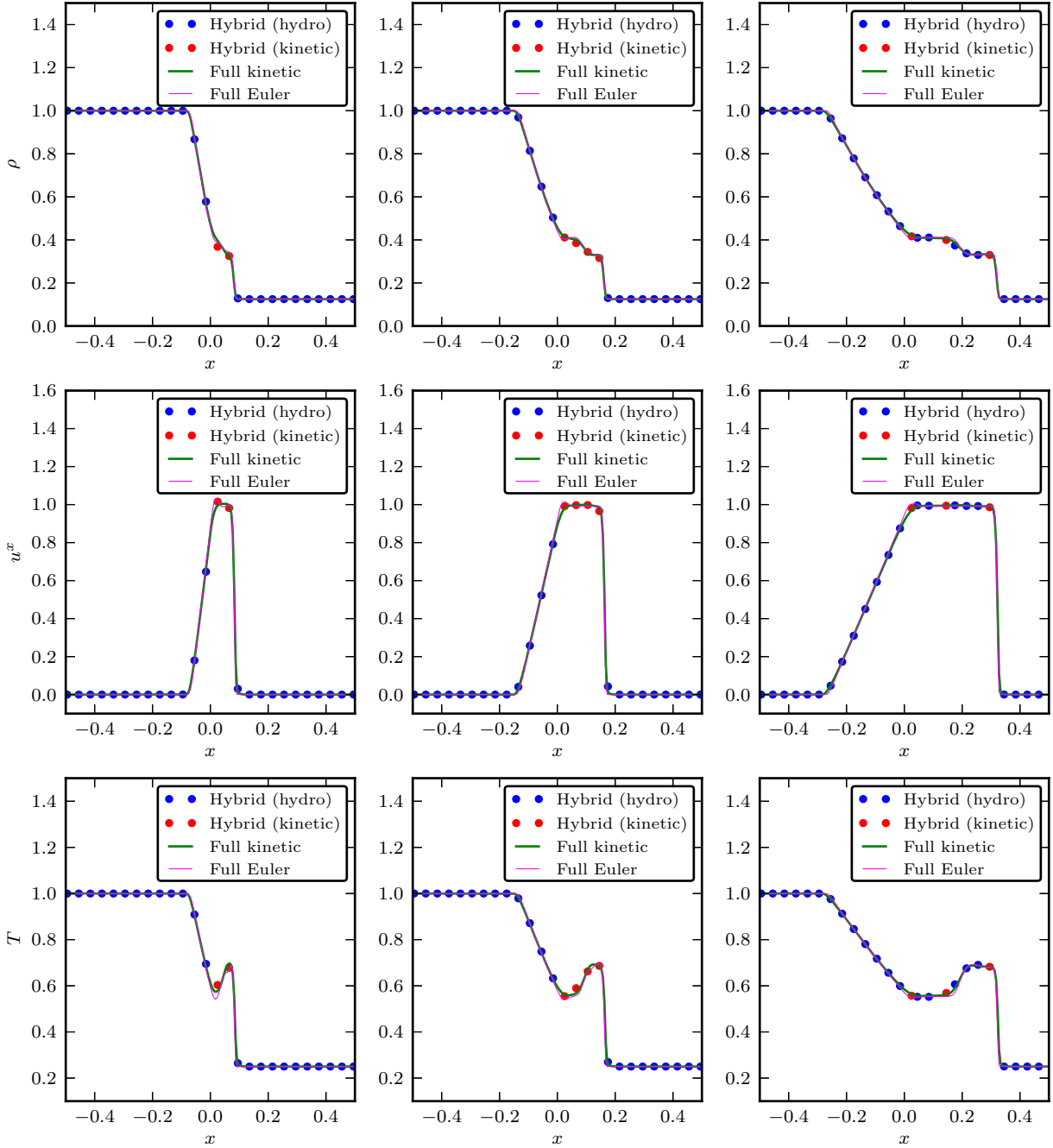}
  \caption{\textbf{Test 1 - Riemann problem with $\ve = 10^{-3}$ :} \emph{Order 0 (Euler)}; Density, mean velocity and temperature at times  $t = 0.05$, $0.10$ and $0.20$.}
  \label{figSodE3-Euler}
\end{figure}

\begin{figure}
  \includegraphics{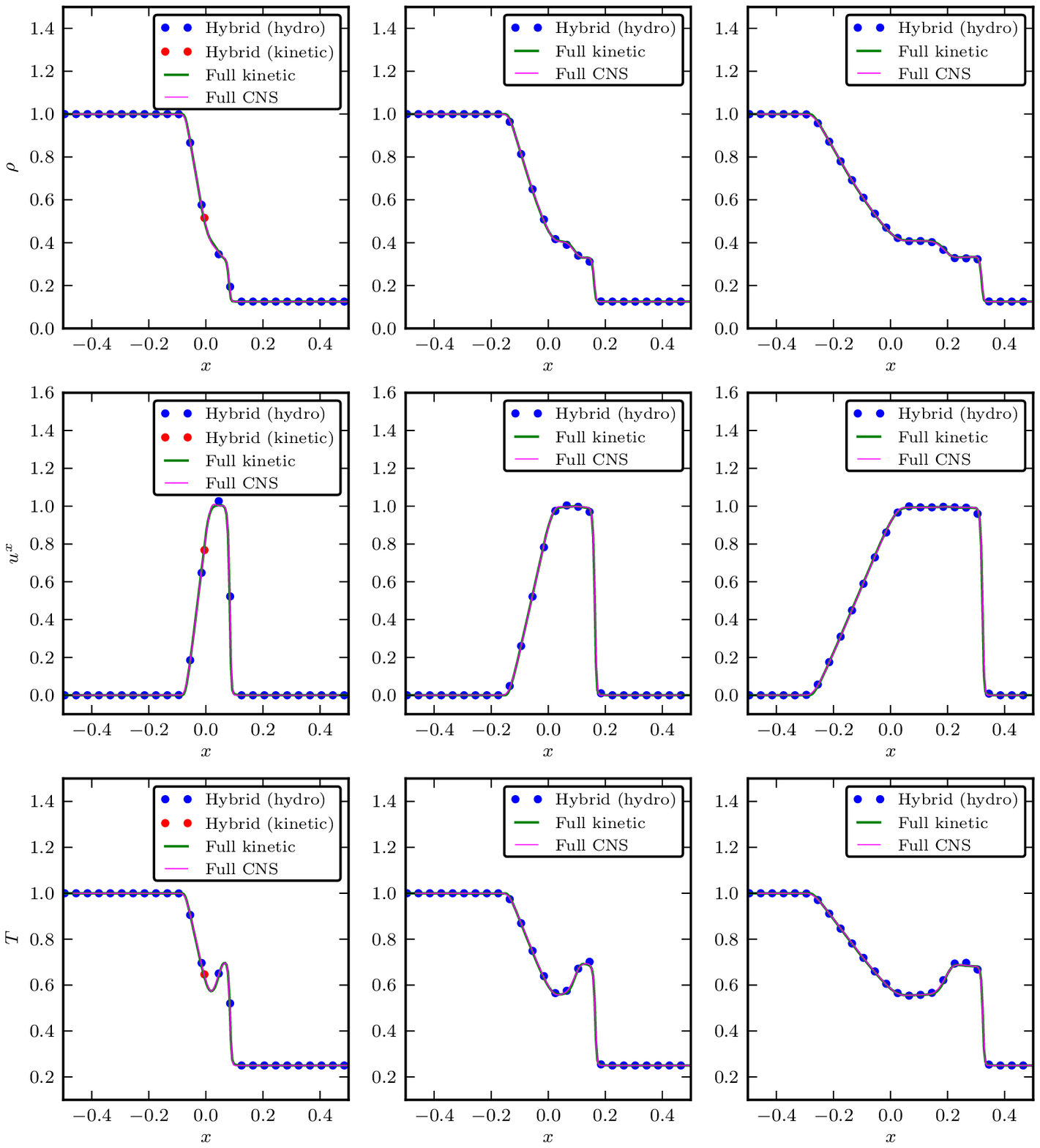}
	\caption{\textbf{Test 1 - Riemann problem with $\ve = 10^{-3}$ :} \emph{Order 1 (CNS)}; Density, mean velocity and temperature at times  $t = 0.05$, $0.10$ and $0.20$.}
	\label{figSodE3-CNS}
\end{figure}

\begin{figure}
  \includegraphics{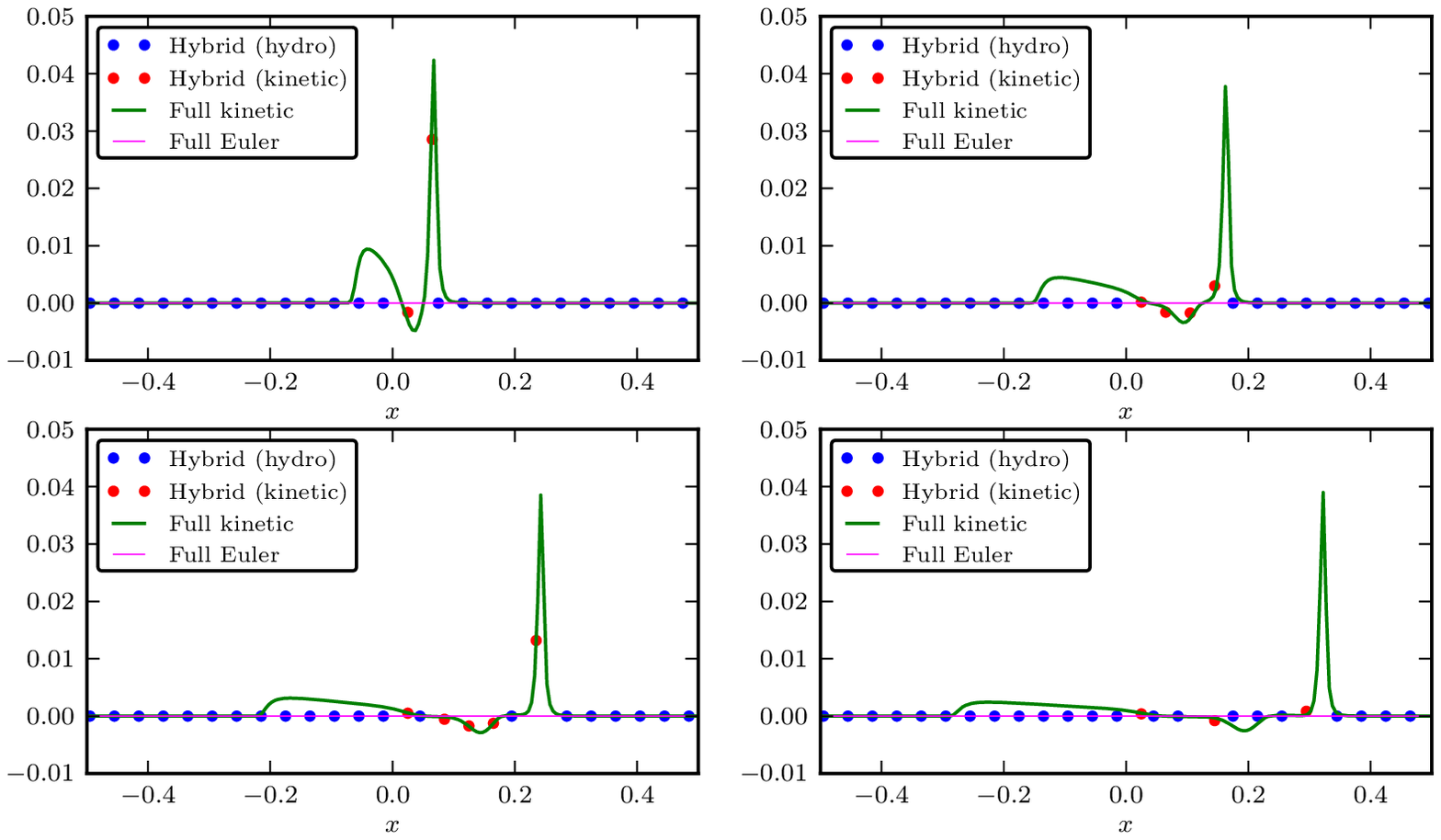}
  \caption{\textbf{Test 1 - Riemann problem with $\ve = 10^{-3}$ :} \emph{Order 0 (Euler)}; heat flux at times  $t = 0.05$, $0.10$, $0.15$ and $0.20$.}
  \label{figSodE3_q-Euler}
\end{figure}

\begin{figure}
  \includegraphics{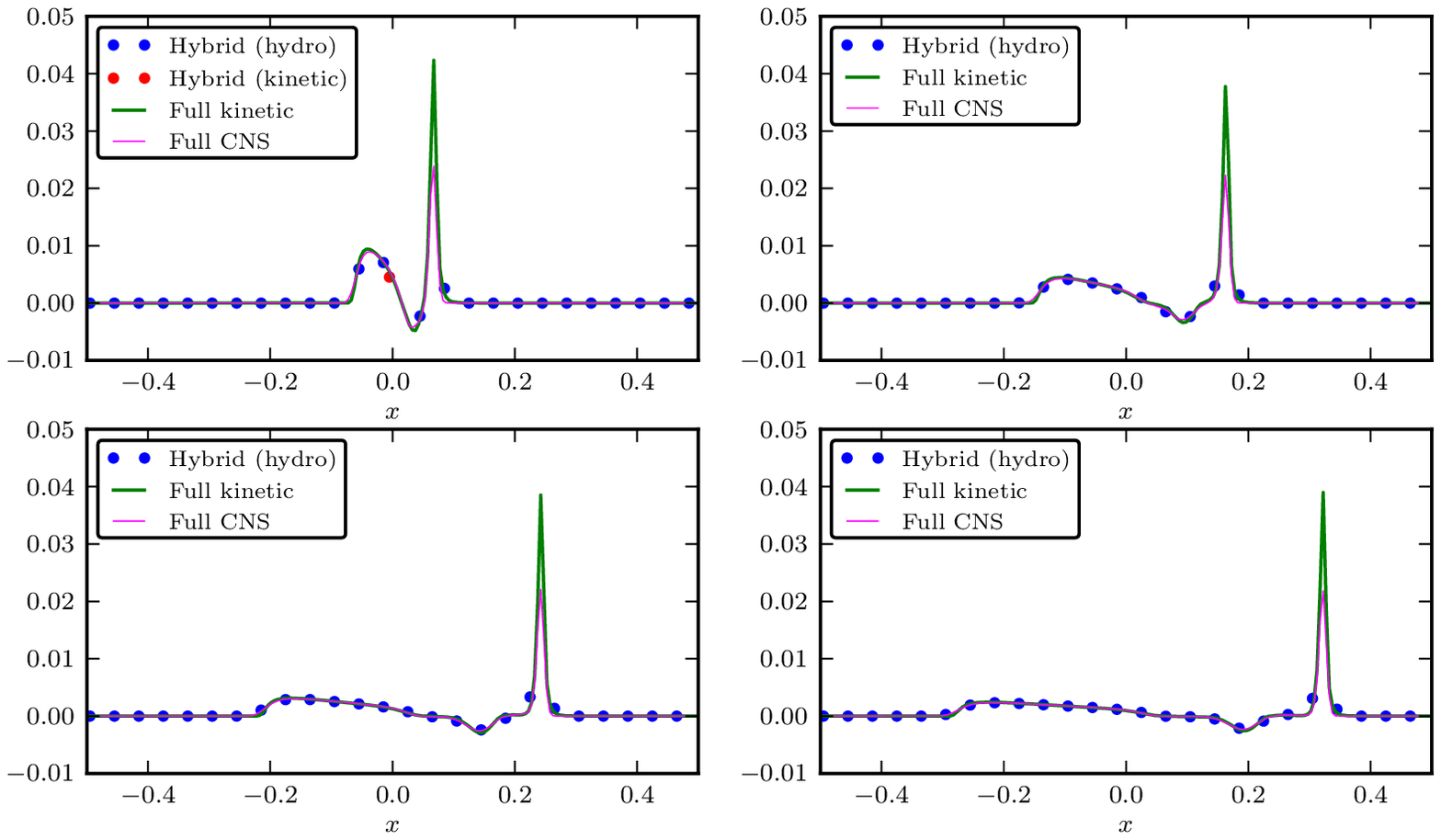}
	\caption{\textbf{Test 1 - Riemann problem with $\ve = 10^{-3}$ :} \emph{Order 1 (CNS)}; heat flux at times  $t = 0.05$, $0.10$, $0.15$ and $0.20$.}
	\label{figSodE3_q-CNS}
\end{figure}

On the other hand, we also give the result of the computations close to the Euler limit ($\ve=10^{-3}$) using $100$ space cells and $32\times 32 \times 32$ cells in velocity for the hybrid method. 
In this case, the solution is very close to the hydrodynamic limit and the kinetic model applies only locally (for instance around a discontinuity where the matrices $\mathcal{V}_{NS}$ and $\mathcal{V}_{Euler}$ differ). Once again, in the Euler case, there is a very good  agreement with the reference solution on the density, mean velocity and temperature   reported in Figures \ref{figSodE3-Euler} and \ref{figSodE3_q-Euler} (although some small error can be seen locally, specially in the heat flux). 
Let us emphasize the the hybrid scheme is perfectly fitted to describe correctly  the time evolution of the heat flux which is zero for the Euler system whereas it fluctuates around zero when the distribution function is not a Maxwellian. For this case, the region where the kinetic models applies is rather small since $\ve \ll 1$ and the hybrid method is particularly efficient. Indeed the computational time is $6$ times faster than  the full kinetic model for the same configuration.
This is even more striking for the first order, CNS case, depicted in Figures \ref{figSodE3-CNS} and \ref{figSodE3_q-CNS}. In this simulation, only one cell is kinetic in short time, and then the whole domain becomes fluid. The error is then negligible (even for the heat flux).
The computational gain becomes huge: the hybrid scheme is $24.5$ times faster than the kinetic one. It becomes almost competitive with the fluid solver.

%
%
%
%

\subsection{Test 2: Blast Wave}
We now consider the case of a blast wave where the initial data is chosen as	
    			      \[ 
			        f^{in}(x,v) = \mathcal{M}_{\rho(x), \bm{u}(x), T(x)}(v), \quad \forall x\in[-0.5,0.5], \quad v\in \RR^3,
			      \]
			      with
			      \[
			        \left (\rho(x), \bm{u}(x), T(x)\right ) = \left \{ 
			          \begin{aligned}
			            & (1,1,0,0,2) && \text{ if } x < -0.3, \\   
			            & (1,0,0,0,0.25) && \text{ if } -0.3 \leq x \leq 0.3, \\   
			            & (1,-1,0,0,2) && \text{ if } x \geq 0.3.
			          \end{aligned}
			        \right.;
			      \]
            Moreover, specular boundary conditions ($\alpha = 1$ in the so-called Maxwellian boundary conditions setting \cite{CIP:94}) are considered in order to study wave reflections.

\begin{figure}
  \includegraphics{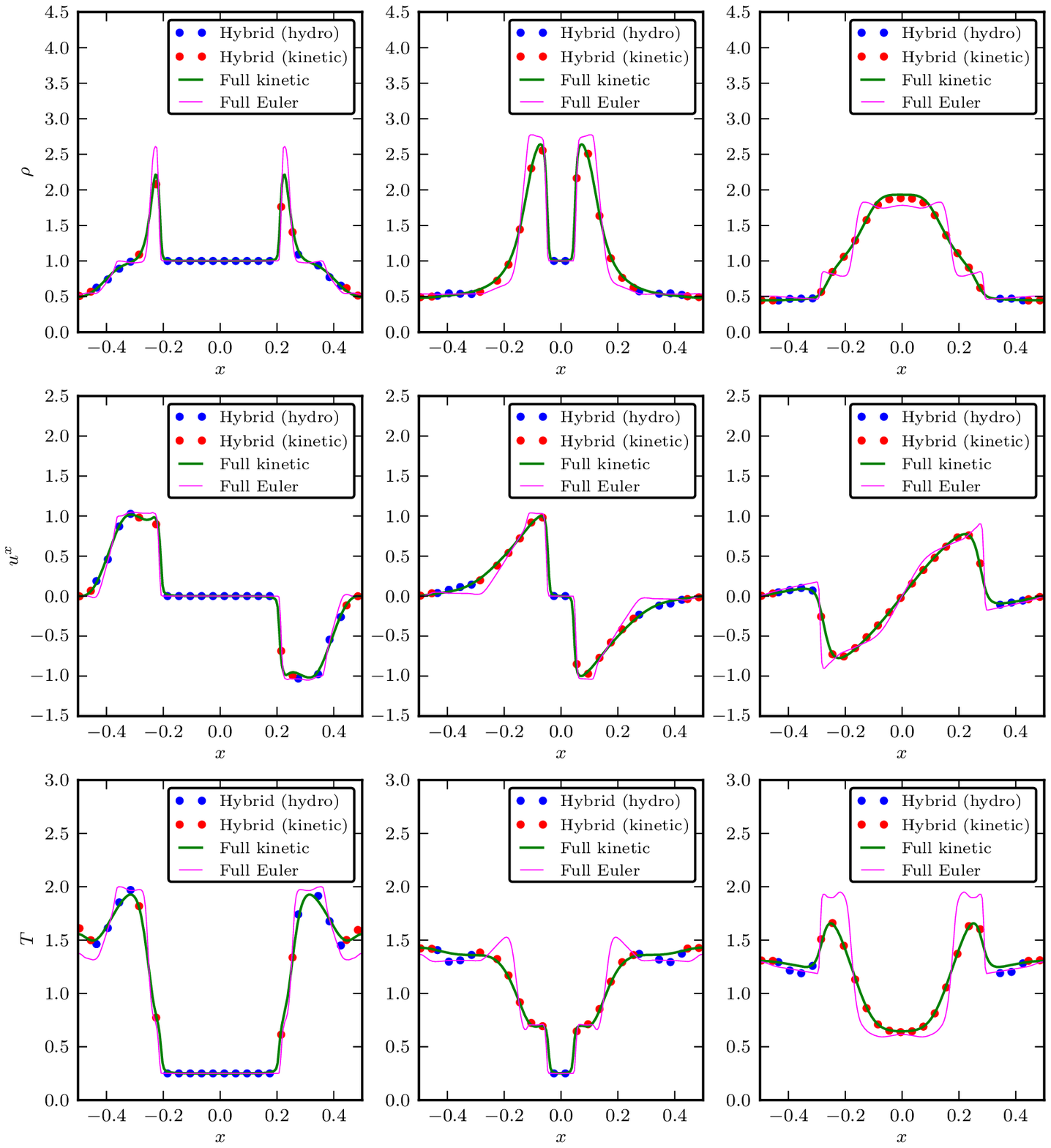}
	\caption{\textbf{Test 2 - Blast wave with $\ve = 10^{-2}$ :} \emph{Order 0 (Euler)}; Density, mean velocity and temperature at times $t = 0.05$, $0.15$ and $0.35$.}
	\label{figBlastE2-Euler}
\end{figure}

\begin{figure}
  \includegraphics{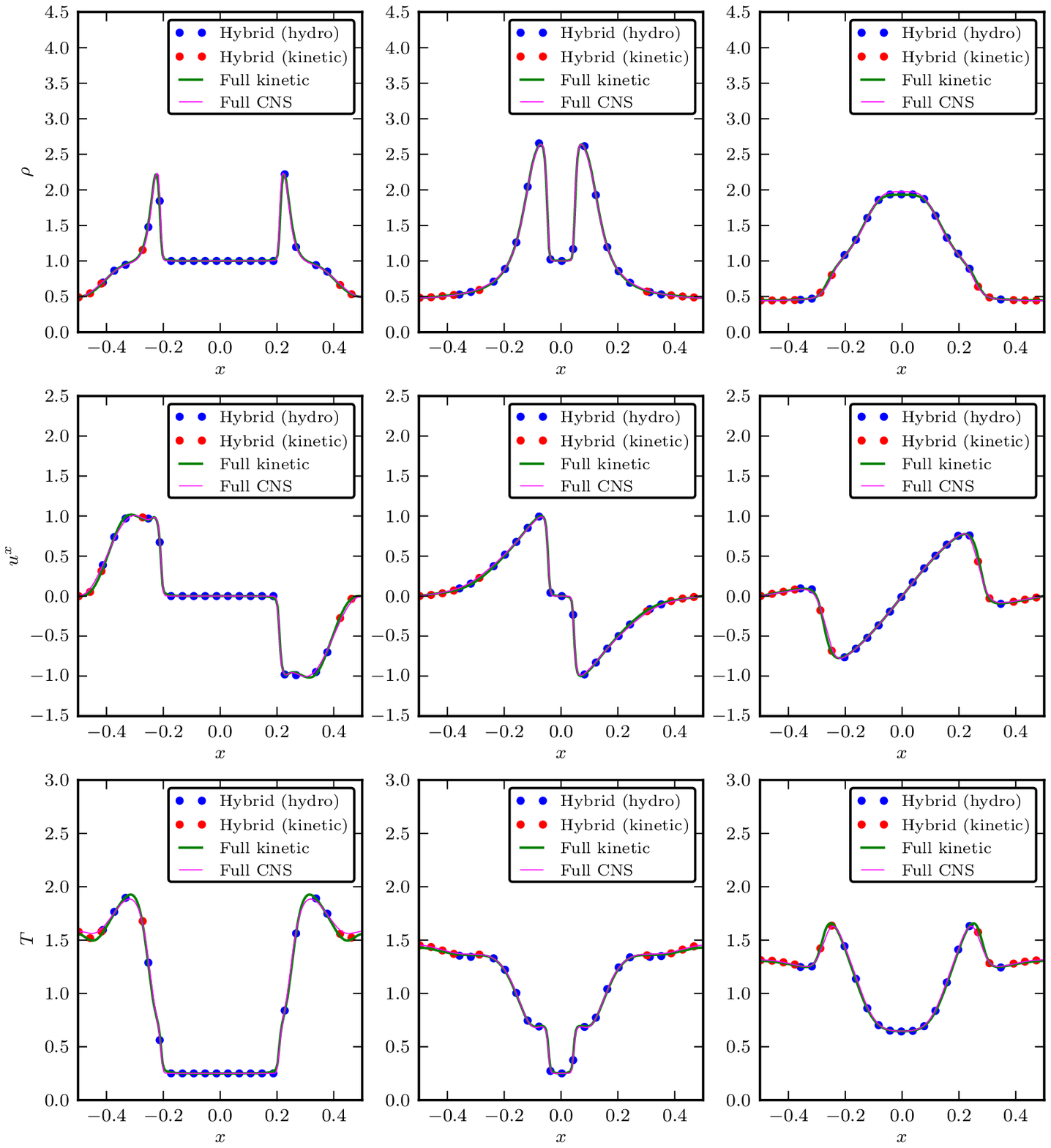}
	\caption{\textbf{Test 2 - Blast wave with $\ve = 10^{-2}$ :} \emph{Order 1 (CNS)}; Density, mean velocity and temperature at times $t = 0.05$, $0.15$ and $0.35$.}
	\label{figBlastE2-CNS}
\end{figure}

\begin{figure}
  \includegraphics{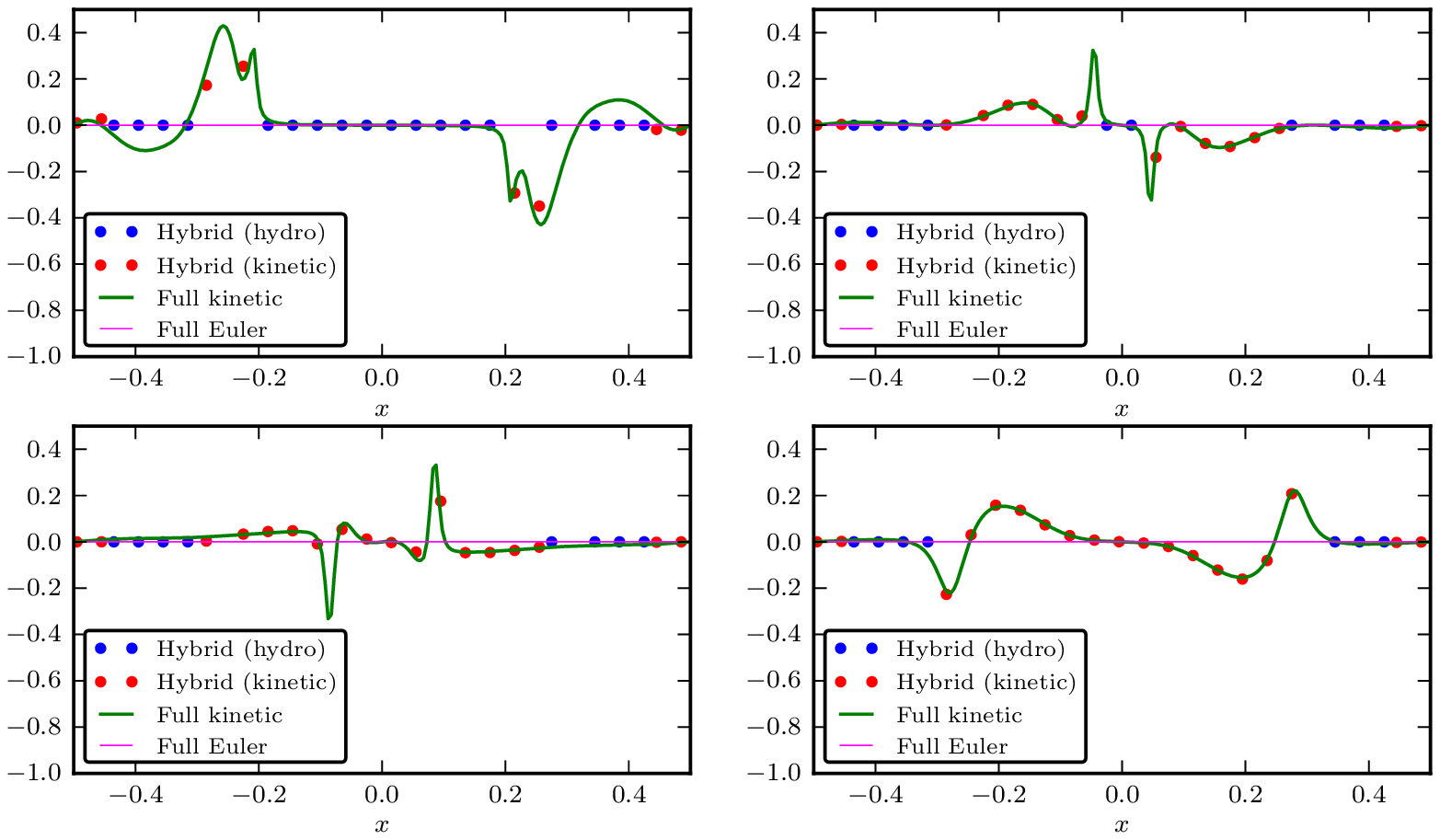}
	\caption{\textbf{Test 2 - Blast wave with $\ve = 10^{-2}$ :} \emph{Order 0 (Euler)}; heat flux at times $t = 0.05$, $0.15$, $0.25$ and $0.35$.}
	\label{figBlastE2-Euler_q}
\end{figure}

\begin{figure}
  \includegraphics{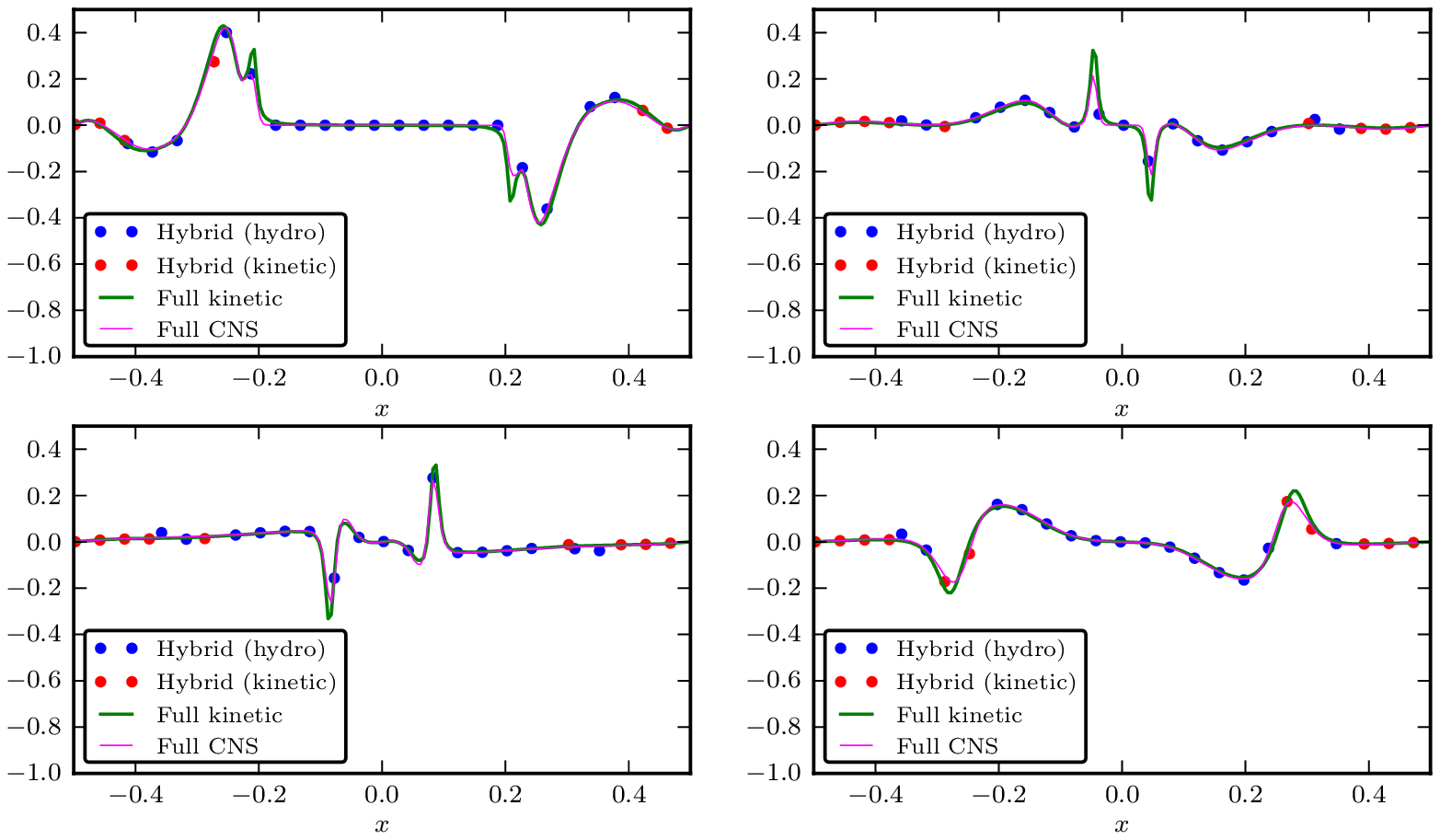}
	\caption{\textbf{Test 2 - Blast wave with $\ve = 10^{-2}$ :} \emph{Order 1 (CNS)}; heat flux at times $t = 0.05$, $0.15$, $0.25$ and $0.35$.}
	\label{figBlastE2-CNS_q}
\end{figure}

            We report the numerical results with $\ve=10^{-2}$ in Figures \ref{figBlastE2-Euler}\,-\, \ref{figBlastE2-CNS_q} at different time with a computational domain in velocity $[-7.5,7.5]^2$. The hybrid scheme is used with $100$ points in $x$ and the size of the velocity grid is $32^3$ points.
One the one hand, the solution obtained with the zeroth order hybrid scheme is compared with ones obtained using a full kinetic model on a fine grid and applying the Euler system. 
We still observe a good agreement between the solution given by the hybrid method and the one obtained with the full kinetic model whereas the  purely macroscopic model does not give accurate results for large time $t \geq 0.1$ (see the Euler case in Figure \ref{figBlastE2-Euler}). The hybrid scheme is quite accurate in the region where the heat flux differs from zero (Figure \ref{figBlastE2-Euler_q}), which confirms the consistency of the criteria described in section~\ref{subEvoCouple}. We nevertheless observe some small discrepancies in large time on these zones.
These errors disappear almost completely when using the first order correction (Figures \ref{figBlastE2-CNS} and \ref{figBlastE2-CNS_q}), even in large time, in particular because the CNS solver is quite close to the kinetic solution. 
   
Concerning the computational time for the same configuration ($N_x=100$ and $N_v=32^3$), the zeroth order hybrid scheme is more efficient than  the kinetic models since the computational time is $2.8$ times faster, and the first order one is $7$ time faster. These improvements are particularly encouraging if we claim to construct an hybrid scheme based on the full Boltzmann operator in $\RR^3_v$ for which the computational complexity is much higher than the BGK operator.

%
%
%
%
%
%
%
%
%
%
%
%

\begin{figure}
  \includegraphics{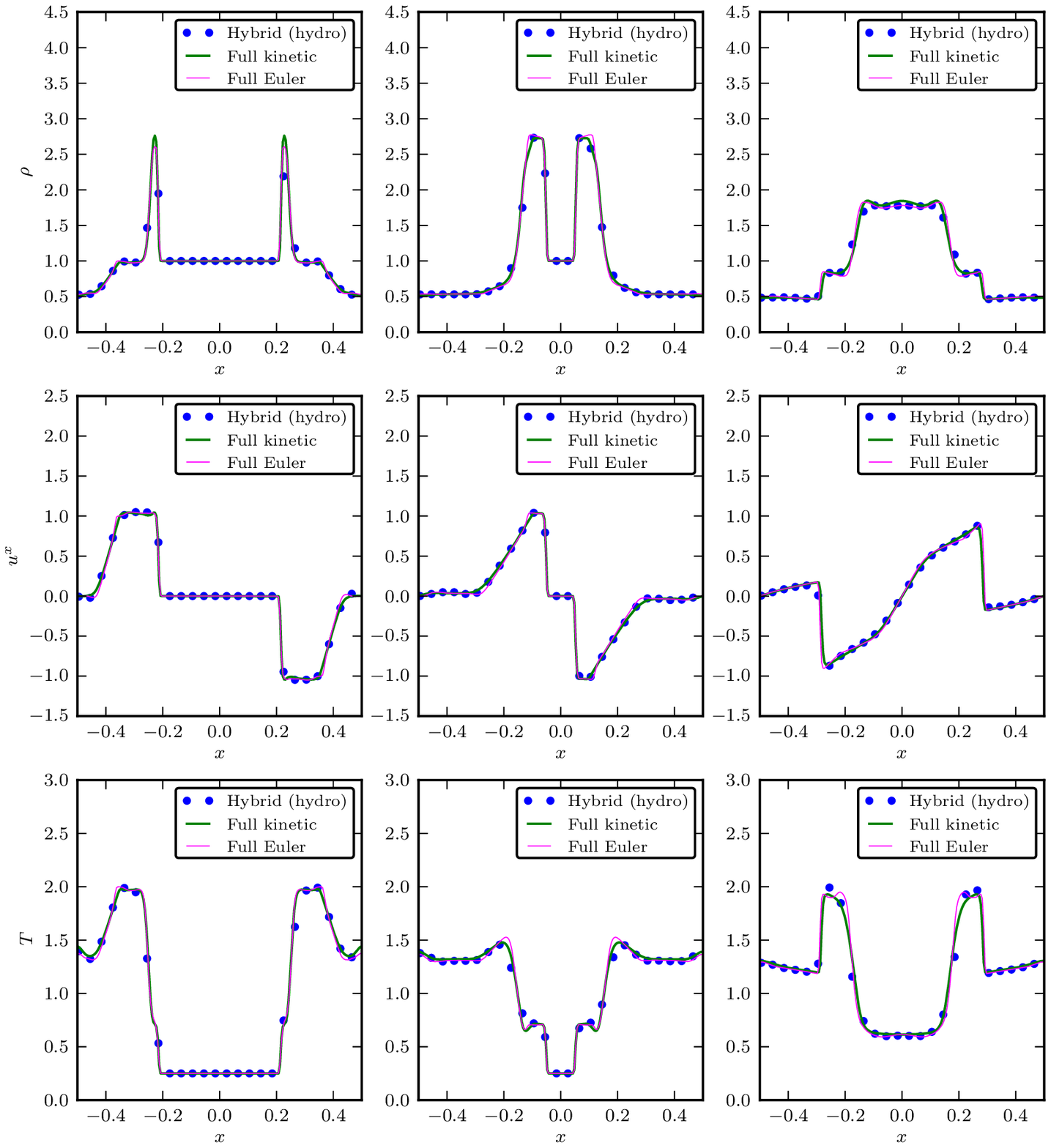}
	\caption{\textbf{Test 2 - Blast wave with $\ve = 10^{-3}$ :} \emph{Order 0 (Euler)}; Density, mean velocity and temperature at times $t = 0.05$, $0.15$ and $0.35$.}
	\label{figBlastE3-Euler}
\end{figure}

\begin{figure}
  \includegraphics{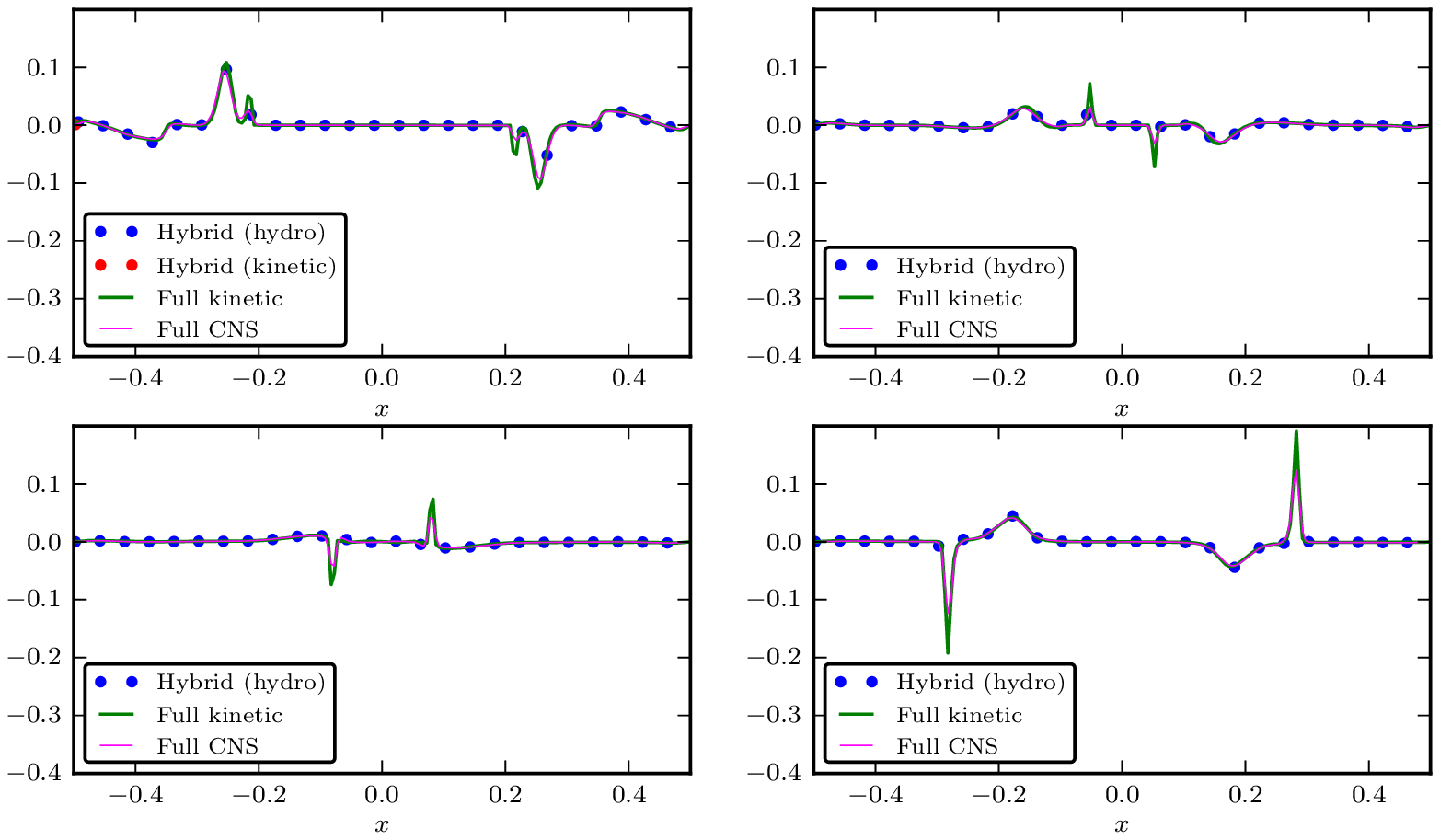}
	\caption{\textbf{Test 2 - Blast wave with $\ve = 10^{-3}$ :} \emph{Order 1 (CNS)}; heat flux at times $t = 0.05$, $0.15$, $0.25$ and $0.35$.}
	\label{figBlastE3-CNS_q}
\end{figure}

On the other hand, we present in Figures \ref{figBlastE3-Euler}\,-\,\ref{figBlastE3-CNS_q} results in the fluid regime $\ve=10^{-3}$. In that case, both zeroth and first order hybrid scheme resolve very accurately the reference solution. The kinetic zones are really small (only the boundary cells are kinetic), and we observe that it is enough for the fluid solvers to achieve the correct result (even if the Euler solver is still far from the reference solution and the heat flux is $0$ everywhere).
The gains in time here are really good, because the Euler solver behave almost like a fluid one (gain of a factor $1250$ compared to the kinetic solver), and the CNS remains competitive, although the parabolic CFL condition has to be applied (gain of a factor $48$).

%
%
%
%
%
\begin{figure}
  \includegraphics{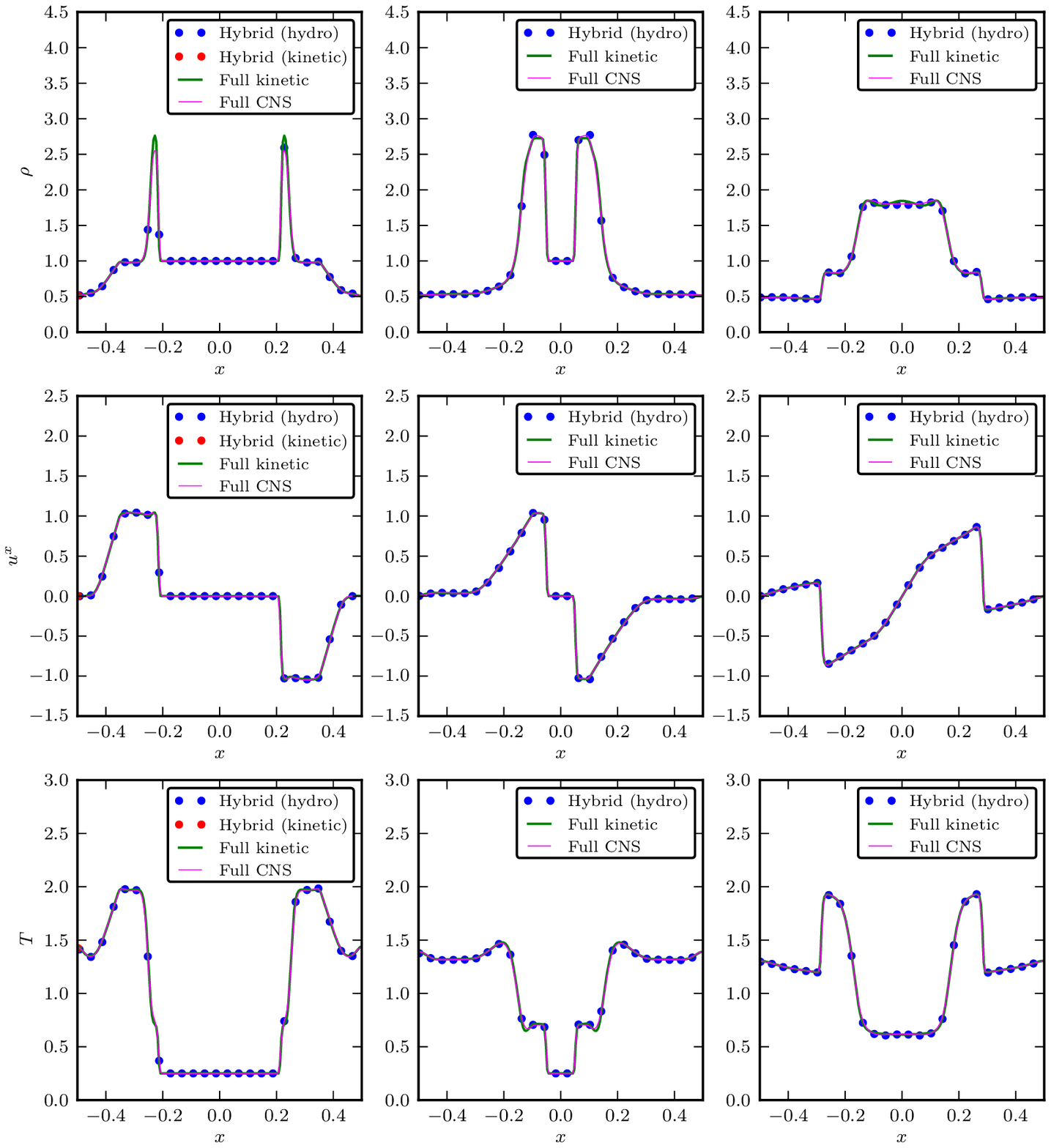}
	\caption{\textbf{Test 2 - Blast wave with $\ve = 10^{-3}$ :} \emph{Order 1 (CNS)}; Density, mean velocity and temperature at times $t = 0.05$, $0.15$ and $0.35$.}
	\label{figBlastE3-CNS}
\end{figure}
%
%
\begin{figure}
  \includegraphics{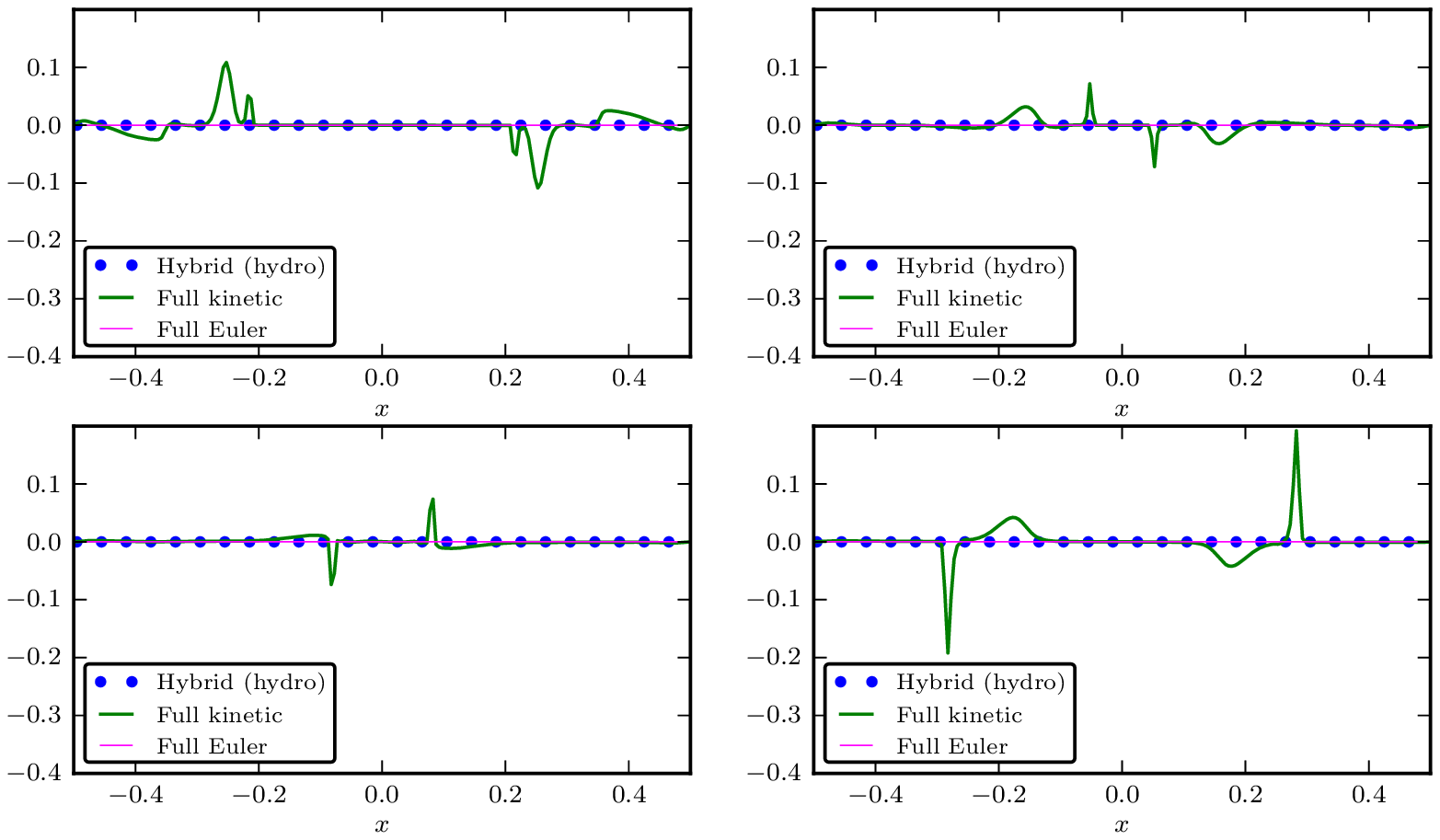}
	\caption{\textbf{Test 2 - Blast wave with $\ve = 10^{-3}$ :} \emph{Order 0 (Euler)}; heat flux at times $t = 0.05$, $0.15$, $0.25$ and $0.35$.}
	\label{figBlastE3-Euler_q}
\end{figure}
%
%
%
%
%

\subsection{Test 3: Far from Equilibrium, Variable Knudsen Number}
This last numerical test deals with the BGK operator where the initial data is far from the thermodynamical equilibrium and when the Knudsen number $\ve$ varies with space. The initial condition is given by 
			      \[ 
			        f^{in}(x,v) = \frac{1}{2}\left ( \mathcal{M}_{\rho(x), \bm{u}(x), T(x)}(v) + \mathcal{M}_{\rho(x), -\bm{u}(x), T(x)}(v) \right ),
			      \]
			      for $x\in[-0.5,0.5], v\in \RR^3$  with
\[
\left (\rho(x), \bm{u}(x), T(x)\right ) \,=\, \left (1 +\frac{1}{2} \sin (\pi x),\, \frac{3}{4},\, 0,\, \frac{5 + 2 \cos( 2\pi x)}{20} \right ).
 \]
Moreover, the Knudsen number $\ve$ varies smoothly from zero to one as 
$$
\ve(x) = 10^{-4} + \frac12\left (\arctan(1+30x) + \arctan(1-30x)\right ).
$$

\begin{figure}
  \includegraphics{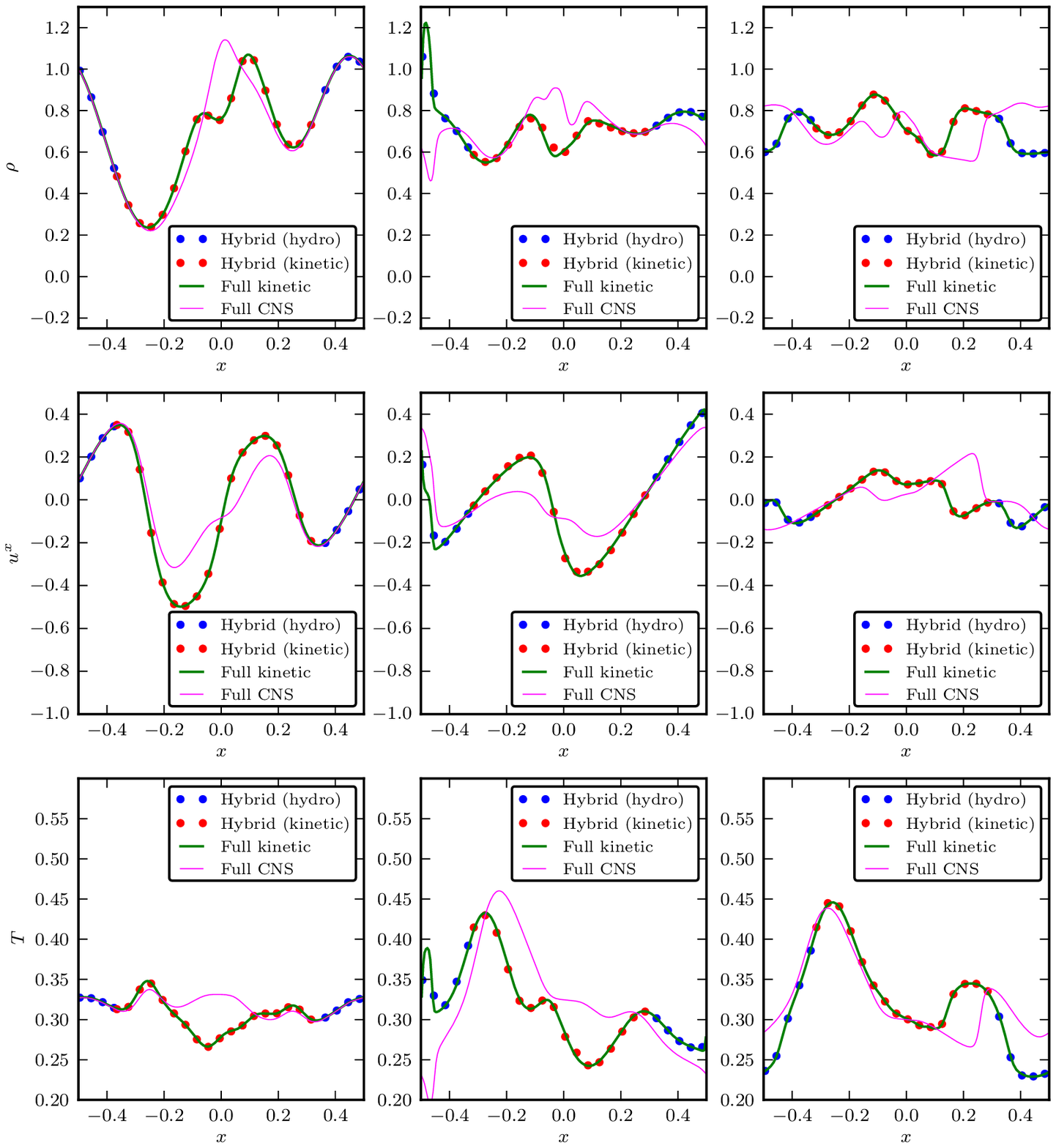}
  \caption{\textbf{Test 3 - Double bump initial data with variable Knudsen number :} \emph{Order 1 (CNS)}; Density, mean velocity and temperature at times $t = 0.10$, $0.50$ and $1.0$.}
  \label{figFar} 
\end{figure}

\begin{figure}
  \includegraphics[scale=1]{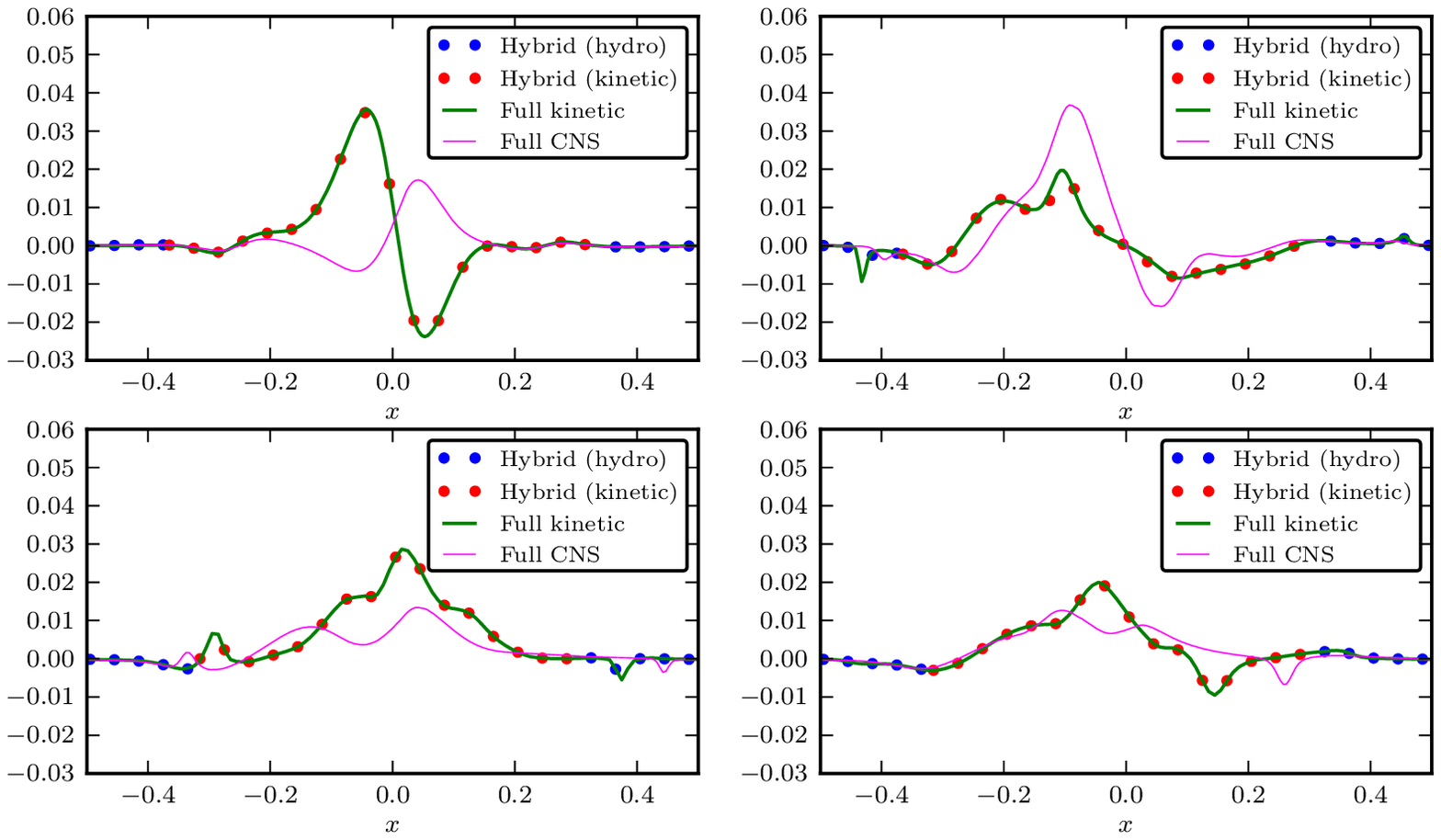}
  \caption{\textbf{Test 3 - Double bump initial data with variable Knudsen number :} \emph{Order 1 (CNS)}; heat flux at times $t= 0.10$, $0.40$, $0.70$ and $1.0$.}
  \label{figFar_q}
\end{figure}

The hybrid scheme is used with $100$ points in $x$ and the size of the velocity grid is $32 \times 32 \times 32$ points on the computational domain $(-8,8)^3$. This numerical test is particularly difficult since the initial data is not at thermodynamical equilibrium and there is no hydrodynamic limit except in the regions where the Knudsen number is small, that is for $|x|\geq 1/3$. 
We compare our numerical solution with the one obtained on a fine grid using the full kinetic model and the one given by solving the compressible Navier-Stokes system. Once again the density, mean velocity and temperature are well described and agree well with the solution corresponding to the kinetic model, even if the fluid model is not correct. 
Indeed, the numerical solution of the fluid equations develops waves propagating in the domain which does not correspond to the solution of the kinetic model (see Figures~\ref{figFar} and \ref{figFar_q}). For such a configuration the hybrid method is $1.9$ times faster than the full kinetic solver.     

%
%
%
%
%
%
%
%

\section{Conclusion}

We propose a simple hierarchy of hybrid method for solving the Boltzmann equation (or analogous kinetic models) in various regimes. This method is based on two criteria. The first one is used to pass from the macroscopic system to the kinetic equation and is strongly inspired by the works of  Levermore, Morokoff and Nadiga in \cite{LevermoreMorokoffNadiga:98}. It is based on a Chapman-Enskog expansion of the distribution. This criterion \eqref{crit:orderK} only depends on macroscopic quantities given by a closure of the kinetic model, and does not require the evaluation of the distribution function.  The second one is used to pass from the kinetic equation to its corresponding hydrodynamical limit and it is based on the comparison of the truncation of the Chapman-Enskog expansion \eqref{devChapEnsk_k} with the exact distribution function with its hydrodynamical equilibrium (\ref{eqKintoFluid1}) and the ratio $\Delta t/\ve$ (\ref{eqKintoFluid2}).   

	\section*{Acknowledgment}
	  The research of the first author (FF) is partially supported by the
European Research Council ERC Starting Grant 2009, project 239983-NuSiKiMo.
	  The research of the second author (TR) was granted by the NSF Grants  \#1008397 and \#1107444 (KI-Net) and ONR grant \#000141210318.
	  TR would like to thanks Dave Levermore for the fruitful discussions and comments about the manuscript.

  \bibliographystyle{acm}
  \bibliography{biblioFR}

\begin{thebibliography}{10}

\bibitem{AlaiaPuppo:2012}
{\sc Alaia, A., and Puppo, G.}
\newblock {A hybrid method for hydrodynamic-kinetic flow, Part II: Coupling of
  hydrodynamic and kinetic models}.
\newblock {\em J. Comput. Phys. 231}, 16 (2012), 5217--5242.

\bibitem{AndriesPerthame:2000}
{\sc Andries, P., {Le Tallec}, P., Perlat, J.-p., and Perthame, B.}
\newblock {The Gaussian-BGK model of Boltzmann equation with small Prandtl
  number}.
\newblock {\em Eur. J. Mech. B Fluids 19}, 6 (Nov. 2000), 813--830.

\bibitem{BerthelinTzavarasVasseur:2009}
{\sc Berthelin, F., Tzavaras, A.~E., and Vasseur, A.}
\newblock From discrete velocity {B}oltzmann equations to gas dynamics before
  shocks.
\newblock {\em J. Stat. Phys. 135}, 1 (2009), 153--173.

\bibitem{Bhatnagar:1954}
{\sc Bhatnagar, P., Gross, E., and Krook, M.}
\newblock {A Model for Collision Processes in Gases. I. Small Amplitude
  Processes in Charged and Neutral One-component Systems}.
\newblock {\em Phys. Rev. 94}, 3 (1954), 511--525.

\bibitem{Boyd:1995}
{\sc Boyd, I.~D., Chen, G., and Candler, G.~V.}
\newblock {Predicting failure of the continuum fluid equations in transitional
  hypersonic flows}.
\newblock {\em Phys. Fluids 7}, 1 (1995), 210.

\bibitem{CIP:94}
{\sc Cercignani, C., Illner, R., and Pulvirenti, M.}
\newblock {\em The Mathematical Theory of Dilute Gases}, vol.~106 of {\em
  Applied Mathematical Sciences}.
\newblock Springer-Verlag, New York, 1994.

\bibitem{DegondDimarco:2012}
{\sc Degond, P., and Dimarco, G.}
\newblock Fluid simulations with localized {B}oltzmann upscaling by direct
  simulation {M}onte-{C}arlo.
\newblock {\em J. Comput. Phys. 231}, 6 (2012), 2414--2437.

\bibitem{DegondDimarcoMieussens:2010}
{\sc Degond, P., Dimarco, G., and Mieussens, L.}
\newblock A multiscale kinetic-fluid solver with dynamic localization of
  kinetic effects.
\newblock {\em J. Comput. Phys. 229}, 13 (2010), 4907--4933.

\bibitem{DimarcoMieussensRispoli:2013}
{\sc Dimarco, G., Mieussens, L., and Rispoli, V.}
\newblock {Asymptotic preserving automatic domain decomposition for the
  Vlasov-Poisson-BGK system with applications to plasmas}.
\newblock Preprint arXiv 1305.1759, 2013.

\bibitem{DimarcoPareschi:2008}
{\sc Dimarco, G., and Pareschi, L.}
\newblock {Hybrid Multiscale Methods II. Kinetic Equations}.
\newblock {\em Multiscale Modeling \& Simulation 6}, 4 (Jan. 2008), 1169--1197.

\bibitem{ellis:1975}
{\sc Ellis, R., and Pinsky, M.}
\newblock {The First and Second Fluid Approximations to the Linearized
  Boltzmann Equation}.
\newblock {\em J. Math. Pures Appl. 54}, 9 (1975), 125--156.

\bibitem{FilbetJin:2010}
{\sc Filbet, F., and Jin, S.}
\newblock {An Asymptotic Preserving Scheme for the ES-BGK Model of the
  Boltzmann Equation}.
\newblock {\em J. Sci. Comput. 46}, 2 (2010), 204--224.

\bibitem{Golse:2005}
{\sc Golse, F.}
\newblock {The Boltzmann equation and its hydrodynamic limits}.
\newblock In {\em Handbook of Differential Equations: Evolutionary Equations
  Vol. 2}, C.~Dafermos and E.~Feireisl, Eds. North-Holland, 2005, pp.~159--303.

\bibitem{Kolobov:2007}
{\sc Kolobov, V., Arslanbekov, R., Aristov, V., a.a. Frolova, and Zabelok, S.}
\newblock {Unified solver for rarefied and continuum flows with adaptive mesh
  and algorithm refinement}.
\newblock {\em J. Comput. Phys. 223}, 2 (May 2007), 589--608.

\bibitem{LevermoreMorokoffNadiga:98}
{\sc Levermore, C.~D., Morokoff, W.~J., and Nadiga, B.~T.}
\newblock Moment realizability and the validity of the {N}avier-{S}tokes
  equations for rarefied gas dynamics.
\newblock {\em Phys. Fluids 10}, 12 (1998), 3214--3226.

\bibitem{NessyahuTadmor:1990}
{\sc Nessyahu, H., and Tadmor, E.}
\newblock {Non-oscillatory central differencing for hyperbolic conservation
  laws}.
\newblock {\em J. Comput. Phys. 87}, 2 (Apr. 1990), 408--463.

\bibitem{Saint-Raymond:2009}
{\sc Saint-Raymond, L.}
\newblock {\em {Hydrodynamic Limits of the Boltzmann Equation}}.
\newblock Springer-Verlag, Berlin, 2009.

\bibitem{Struchtrup:2005}
{\sc Struchtrup, H.}
\newblock {\em {Macroscopic Transport Equations for Rarefied Gas Flows}}.
\newblock Springer-Verlag, Berlin, 2005.

\bibitem{Tiwari:98}
{\sc Tiwari, S.}
\newblock Coupling of the {B}oltzmann and {E}uler equations with automatic
  domain decomposition.
\newblock {\em J. Comput. Phys. 144}, 2 (1998), 710--726.

\bibitem{Tiwari:2000}
{\sc Tiwari, S.}
\newblock {Application of moment realizability criteria for the coupling of the
  Boltzmann and Euler equations}.
\newblock {\em Transport Theory Statist. Phys. 29}, 7 (2000), 759--783.

\bibitem{TiwariKlarHardt:2009}
{\sc Tiwari, S., Klar, A., and Hardt, S.}
\newblock A particle-particle hybrid method for kinetic and continuum
  equations.
\newblock {\em J. Comput. Phys. 228}, 18 (2009), 7109--7124.

\bibitem{TiwariKlarHardt:2012}
{\sc Tiwari, S., Klar, A., and Hardt, S.}
\newblock {Simulations of micro channel gas flows with domain decomposition
  technique for kinetic and fluid dynamics equations}.
\newblock In {\em 21st International Conference on Domain Decomposition
  Methods\/} (2012), pp.~197--206.

\end{thebibliography}

\end{document}